\documentclass[11pt]{article}

       \font\tenmsb=msbm10
       \font\sevenmsb=msbm7
       \font\fivemsb=msbm5
       \catcode`\@=11
       \ifx\amstexloaded@\relax\catcode`\@=\active
       \endinput\else\let\amstexloaded@\relax\fi
       \def\spaces@{\space\space\space\space\space}
       \def\spaces@@{\spaces@\spaces@\spaces@\spaces@\spaces@}
       \def\space@.{\futurelet\space@\relax}
       \space@. %
       \def\Err@#1{\errhelp\defaulthelp@\errmessage{AmS-TeX error: #1}}
       \def\relaxnext@{\let\next\relax}
       \def\accentfam@{7}
       
       \def\noaccents@{\def\accentfam@{0}}
       \def\Cal{\relaxnext@\ifmmode\let\next\Cal@\else
       \def\next{\Err@{Use \string\Cal\space only in math mode}}\fi\next}
       \def\Cal@#1{{\Cal@@{#1}}}
       \def\Cal@@#1{\noaccents@\fam\tw@#1}

       \def\Bbb{\relaxnext@\ifmmode\let\next\Bbb@\else
       \def\next{\Err@{Use \string\Bbb\space only in math mode}}\fi\next}
       
       \def\Bbb@#1{{\Bbb@@{#1}}}
       \def\Bbb@@#1{\noaccents@\fam\msbfam#1}
       \newfam\msbfam
       \textfont\msbfam=\tenmsb
       \scriptfont\msbfam=\sevenmsb
       \scriptscriptfont\msbfam=\fivemsb
\usepackage{hyperref}
\usepackage{tikz}
\usepackage{pgflibraryarrows}
\usepackage{pgflibrarysnakes}
\usepackage{amsmath,amssymb}

\def\E{{\Bbb E}}
\def\N{{\Bbb N}}
\def\Z{{\Bbb Z}}
\def\R{{\Bbb R}}
\def\T{{\Bbb T}}
\def\C{{\Bbb C}}
\def\Q{{\Bbb Q}}

\newtheorem{thm}{Theorem}
\newtheorem{thm*}{Theorem}
\newtheorem{lemma}{Lemma}[section]
\newtheorem{prop}{Proposition}

\newtheorem{cor}{Corollary}
\newtheorem{claim}{Claim}
\newtheorem{rmk}{Remark}[section]

\newtheorem{notation*}{Notation}

\@addtoreset{equation}{section}

\newcommand{\qed}{\nolinebreak\hfill\rule{2mm}{2mm}
\par\medbreak}
\newcommand{\proof}{\par\medbreak\it Proof: \rm}

\newcommand{\beq}{\begin{equation} }
\newcommand{\eeq}{\end{equation} }

\begin{document}
\title{Quantitative continuity of singular continuous
  spectral measures and arithmetic
criteria 
for quasiperiodic
Schr\"odinger operators.
}
\author{Svetlana Jitomirskaya, Shiwen Zhang}
\date{}
\maketitle
\begin{abstract}
We introduce a notion of $\beta$-almost periodicity and prove
quantitative lower spectral/quantum dynamical bounds for general
bounded $\beta$-almost periodic potentials. Applications include a sharp arithmetic criterion of full spectral dimensionality for analytic quasiperiodic
Schr\"odinger operators in the positive Lyapunov exponent regime  and  arithmetic
criteria  for families with zero Lyapunov exponents, with applications to
Sturmian potentials and the
critical almost Mathieu operator.
\end{abstract}

\tableofcontents
\newpage

\section{Introduction}

Singular continuous spectral measures of Schrodinger operators, usually defined by what they are
not, 
are still not very well understood. The aim of direct spectral theory
is to obtain properties of spectral measures/spectra
and associated quantum dynamics based on the properties of the
potential. In the context of 1D operators this is most often done via
the study of solutions/transfer matrices/dynamics
of transfer-matrix cocycles. Indeed, there are many beautiful results
linking the latter to either  dimensional properties of spectral
measures (going back to \cite{jl1}) or directly to  quantum dynamics
(e.g. \cite{kkl,dt}).  There is also a long thread of results relating
dimensional properties of spectral measures to quantum dynamics
(e.g. \cite{bsb,bgt} and references therein) as well as results
connecting spectral/dynamical properties to some further aspects
(e.g. \cite{last93,bls}). Many of those have been used to obtain
dimensional/quantum dynamical results (sometimes sharp) for several
concrete families (e.g. \cite{dgy}). However, there were no results
directly linking easily formulated properties of the potential to
dimensional/quantum dynamical results, other than for specific
families or a few that ensure
either the mere singularity or continuity of spectral measures (and their
immediate consequences). In particular, we don't know of any
quantitative results of this type.

In this paper we prove the first such result. Consider Schr\"odinger operator on $\l^2(\Z)$ given by
 \begin{equation}\label{Schro-op-gen}
    (H u)_n= u_{n+1}+u_{n-1} + V(n)u_n
\end{equation}
For $\beta>0,$ we  say a real sequence $\{V(n)\}_{n\in\Z}$ has $\beta$-repetitions if there is a sequence of positive integers $q_n\to\infty$ such that
\begin{equation}\label{betaRP}
    \max_{1\le j \le q_n}|V(j)-V(j\pm q_n)|\le e^{-\beta q_n}
\end{equation}
 We will say that $\{V(n)\}_{n\in\Z}$ has $\infty$-repetitions if
(\ref{betaRP}) holds for any $\beta>0$.
For $\beta<\infty,$ we will say that
$\{V(n)\}_{n\in\Z}$ is {\it $\beta$-almost periodic,} if, for some
$\epsilon>0,$  $V(\cdot+kq_n)$ satisfies (\ref{betaRP}) for any $|k|\le e^{\epsilon\beta q_n}/q_n$, i.e., 
\begin{equation}\label{beta-exp-peri}
    \max_{1\le j \le q_n, |k|\le e^{\epsilon\beta  q_n}/q_n}|V(j+kq_n)-V(j+(k\pm1)
q_n)|\le e^{-\beta q_n}
\end{equation}
for any $n.$  
We will say that $\{V(n)\}_{n\in\Z}$ is {\it $\infty$-almost
  periodic,} if it is $\beta$-almost periodic for any $\beta<\infty.$
We note that $\beta$ and even $\infty$-almost periodicity does not
imply almost periodicity in the usual sense. In particular, it is
easily seen that there is an explicit set of generic skew shift potentials that
satisfy this condition.

We will prove

\begin{thm}\label{00}Let $H$ be given by  (\ref{Schro-op-gen})
and  $V$ is bounded and  $\beta$-almost periodic. 
Then, for  an explicit $C=C(\epsilon,V)>0$, 
for any
\begin{equation}
    \gamma<{1-C/\beta}
\end{equation}
the spectral measure  is $\gamma$-spectral continuous. 
\end{thm}
For the definition of spectral continuity (a property that also implies packing continuity
and thus lower bounds on quantum
dynamics) see Section \ref{main}. We formulate a 
more precise (specifying the dependence
of $C$ on $\epsilon,V$) version
in Theorem
\ref{thmConti}.

Our result can be viewed as a quantitative version
simultaneously of two well known statements
\begin{itemize}
\item {\it  Periodicity implies absolute continuity}. Indeed, we prove that a
  quantitative weakening ($\beta$-almost periodicity) implies
  quantitative continuity of the (fractal) spectral measure.
\item{\it Gordon condition (a single/double almost repetition) implies
    continuity of the spectral measure}. Indeed, we prove that a quantitative
 strengthening (multiple almost repetitions) implies  quantitative continuity
 of the spectral measure.
\end{itemize}

Potentials with $\infty$-repetitions are known in the literature as Gordon
potentials \footnote{While $\infty$-repetitions are
      usually used in the definition of Gordon potentials, typically $\beta$-repetitions for
      sufficiently large $\beta$ are enough for the applications}. This property has been used  fruitfully in
the spectral theory in various situations, see reviews \cite{damgorrev,damrev} and references therein. In many cases those
potentials were automatically $\beta$ or even $\infty$- almost
periodic, so satisfied almost repetitions over sufficiently many
periods.  However, even in such cases, what all
those papers used was the strength of the approximation over one-two
(almost) periods based on Gordon Lemma type arguments. Our main technical accomplishment here is that we
find a new algebraic argument and develop technology that allows to
obtain quantitative corollaries from
the fact that the approximation stays strong over many periods, thus
exploring this feature analytically for the first time.

Lower bounds on spectral dimension lead to lower bounds on
packing dimension, thus also for the packing/upper box counting dimensions of the spectrum
as a set and for the upper rate of quantum dynamics. Therefore we obtain
corresponding non-trivial results for all above quantities.

It is clear that our general result only goes in one direction, as
even absolute continuity of the spectral measures does not imply
$\beta$-almost periodicity for $\beta>0.$

However, in the important context of analytic quasiperiodic operators this leads to a sharp
if-and-only-if result.

Let $H=H_{\theta,\alpha,V}$ be a Schr\"odinger operator on $\l^2(\Z)$ given by
\begin{equation}\label{Schro-op}
    (H u)_n= u_{n+1}+u_{n-1} + V
(\theta+n\alpha)u_n, \ \  n\in\Z,\ \ \theta\in\T
\end{equation}
where $V$ is the potential, $\alpha\in\R\backslash\Q$ is the frequency
and $\theta\in\T$ is the phase. Let $\mu=\mu_{\theta,\alpha}$ be the
spectral measure associated with vectors $\delta_0,\delta_1\in
l^2(\Z)$ in the usual sense.

Given $\alpha\in (0,1)$, let ${p_n}/{q_n}$ be the continued fraction approximants to $\alpha$. Define
\begin{equation}\label{beta}
    \beta(\alpha):=\limsup_n\frac{\log q_{n+1}}{q_n}\in[0,\infty].
\end{equation}

Let $S:=\{E\in \sigma(H): L(E)>0\},$ where $\sigma(H)$ is the spectrum
of $H$ and $L(E)$ is the Lyapunov exponent, be the set of
supercritical energies (or, equivalently, the set of $E$ such that the
corresponding transfer-matrix cocycle is non-uniformly
hyperbolic). $S$ depends on $\alpha$ and $V$ but not on $\theta.$

Our main application is
\begin{thm}\label{0} For any analytic $V$ and any $\theta,$ the spectral
  measure $\mu$ restricted to $S$ is of full spectral dimension if and
only if $\beta(\alpha)=\infty.$
\end{thm}

 Full spectral
dimensionality is defined through the boundary behavior of Borel
transform of the spectral measure (see details in Section \ref{main}). It implies a range of properties,
in particular, maximal packing dimension and quasiballistic quantum
dynamics. Thus our criterion links in a sharp way a purely analytic
property of the spectral measure to arithmetic property of the
frequency. The result is local (so works for any subset of the
supercritical set, see Theorem \ref{thmmain} for more detail) and
quantitative (so we obtain separately quantitative spectral
singularity and spectral continuity statements for every finite value of
$\beta,$ see Theorems \ref{thmSingular} and \ref{thmConti}).

The study of one-dimensional one-frequency quasiperiodic operators with general analytic
potentials has seen remarkable advances in the last two decades, from
the Eliasson's  KAM point spectrum proof for the general class  \cite{eli}, to
Bourgain-Goldstein's non-perturbative method \cite{bg}, to Avila's
global theory \cite{a-glob}. In particular, many results have been obtained
in the regime of positive Lyapunov exponents (dubbed
 supercritical in \cite{a-glob}). They can be divided into two
 classes
\begin{itemize}
\item Those that hold for all frequencies (e.g. \cite{jl2,bj02,dt,jmx,jmv1,jmv2})
\item Those that have arithmetic (small denominator type) obstructions
 preventing their holding for all frequencies thus requiring a Diophantine type
  condition (e.g. \cite{bg,
gs2,dgsv}) \footnote{Not all  results can currently be classified this way, most notable
example being the Cantor structure of the spectrum \cite{gs3}, currently proved
for a non-arithmetically-defined full measure set of frequencies, while
the statement has no known arithmetic obstructions. Theoretically
there may also be results such as \cite{ajten} which formally should
belong to the first group but the proof requires argument that highly
depends on the arithmetics, so they must be in the second group, in
spirit. In some sense \cite{bj02} is a result of this type.}
\end{itemize}
Results of the first kind often (but not always \cite{bj02,wy}) do
not require analyticity and hold in higher generality. Results of
the second kind describe phenomena where there is a transition in the
arithmetics of the frequency, thus an extremely interesting question is to
determine where does this transition happen and to understand the neighborhood
of the transition. However, even though some improvements on the
frequency range of some results above have been obtained
(e.g. \cite{yz}), most
existing proofs often require a removal of a
non-arithmetically defined measure zero set of frequencies, thus
cannot be expected to work up to the transition. There have been
remarkable recent advances in obtaining complete arithmetic criteria in
presence of transitions \cite{ayz,jliu1,jliu2} or non-transitions
\cite{ajten} for explicit  popular Hamiltonians: almost Mathieu operator and
Maryland model, but there have been no
such results that work for large families of
potentials. Theorem \ref{0} is  the first theorem of this kind.

A natural way to distinguish between different singular continuous
spectral measures
is by their Hausdorf dimension. However Hausdorff dimension is a poor tool for
characterizing the singular continuous spectral measures arising in
the regime of positive Lyapunov exponents, as it is always equal to
zero (for a.e. phase for any ergodic case \cite{sim}, and for every phase for
one frequency analytic potentials \cite{jl2}\footnote{The result of \cite{jl2}  is
  formulated for trigonometric polynomial $v$. However it extends to
  the analytic case - and more - by the method of \cite{jmv2}.}). Similarly, the lower transport
exponent is always zero for piecewise Lipshitz potentials
\cite{dt,jmv2}. Thus those
two quantities don't even distinguish between pure point and singular
continuous situations. In contrast, our quantitative version of
Theorem \ref{0}, contained in  Theorems \ref{thmSingular} and \ref{thmConti}, shows that spectral dimension is
a good tool to finely distinguish between different kinds of singular
continuous spectra appearing in the supercritical regime for analytic potentials.

The continuity part of Theorem \ref{0} is robust and only requires
some regularity of $V$.
Besides the mentioned
criterion,
Theorem \ref{00}  allows us to obtain new
results for other popular models, such as the critical almost Mathieu
operator, Sturmian potentials, and others.

Indeed, 
our lower bounds are effective for $\beta>C\sup_{E\in\sigma(H)}L(E)$ where $L(E)$
is the Lyapunov exponent (see Theorem
\ref{thmConti}) thus the
range of $\beta$  is increased  for smaller Lyapunov exponents, and in
particular, we obtain full spectral dimensionality (and therefore
quasiballistic motion) as long as $\beta(\alpha)>0,$ when Lyapunov
exponents are zero on the spectrum. This applies, in particular, to
Sturmian potentials and the critical almost Mathieu operator.

As an
example, setting $S_0=\{E: L(E)=0\}$ 
we have
\begin{thm}\label{lip}  For Lipshitz $V$, the quantum dynamics is quasiballistic
\begin{enumerate}\item for any $\beta(\alpha)>0,$ if $S_0\not=\emptyset$
\item for $\beta(\alpha)=\infty,$ otherwise
\end{enumerate}
\end{thm}

A similar statement also holds for full spectral dimensionality or
packing/box counting dimension one. The Lipshitz condition can be relaxed to
piecewise Lipshitz (or even H\"older), leading to part 1 also holding for Sturmian
potentials. This in turn leads to first explicit examples of operators
whose integrated density of state has different Hausdorff and packing dimensions, within both the critical almost Mathieu and
Sturmian families.

The fact that quantum motion can be quasiballistic for highly
Liouville frequencies was first realized by Last \cite{l1} who
proved that almost Mathieu operator with an appropriate (constructed
step by step)  Liouville frequency is quasiballistic. Quasiballistic
property is a $G_\delta$ in any regular (a-la Simon's
Wonderland theorem \cite{sim1}) space \cite{gkt,co}, thus this was known
for (unspecified)  topologically generic frequencies.  Here we show
a precise arithmetic condition on $\alpha$ depending on whether or not Lyapunov
exponent vanishes somewhere on the spectrum.
Thus, in the regime of positive Lyapunov exponents, the quantum motion is very interesting, with dynamics almost
bounded along some scales \cite{jmv2} (this property is sometimes called quasilocalization) and almost ballistic along others.
For finite values of $\beta(\alpha)$ in this regime  our result also yields power-law quantum
dynamics along certain scales while bounded along others.


\subsection{Main application}\label{main}
Fractal properties of Borel measures on $\mathbb{R}$ are linked to the
boundary behavior of their Borel transforms \cite{djls}.  Let
\begin{equation}\label{M}
    M(E+ i\varepsilon)=\int\frac{{\rm d}\mu(E')}{E'-(E+ i\varepsilon)}
\end{equation}
be the Borel transform of measure $\mu$. Fix $0<\gamma<1$. If for $\mu$ a.e. $E$,
\begin{equation}\label{gammaC}
    \liminf_{\varepsilon\downarrow0}\varepsilon^{1-\gamma}|M(E+ i\varepsilon)|<\infty,
\end{equation}
we say measure $\mu$ is (upper) $\gamma$-spectral continuous. Note
that spectral continuity (and singularity) captures the $\liminf$
power law behavior of $M(E+ i\varepsilon)$, while the corresponding
$\limsup$ behavior is linked to the Hausdorff dimension \cite{djls}.
Define the (upper) spectral dimension
 of $\mu$ to be
\begin{equation}\label{sd}
    s(\mu)=\sup\big\{\gamma\in(0,1):\ \mu\ \textrm{is}\ \gamma\textrm{-spectral continuous}\big\}.
\end{equation}
For a Borel subset $S\subset \mathbb{R},$  let $\mu_{S}$ be the
restriction of $\mu$ on $S$. A reformulation of Theorem \ref{0} is
\begin{thm}
\label{thmmain}
Suppose $V$ is real analytic and $L(E)>0$ for every $E$ in some Borel set $S\subset\R$.
Then for any $\theta\in\T$,
$s(\mu_{S})=1$ if and only if $\beta(\alpha)=+\infty$.
\end{thm}
\begin{rmk}
If for $\mu$ a.e. $E$,
\begin{equation}\label{gammaS}
    \liminf_{\varepsilon\downarrow0}\varepsilon^{1-\gamma}|M(E+ i\varepsilon)|=+\infty,
\end{equation}
we say measure $\mu$ is (upper) $\gamma$-spectral singular. We can also consider
\begin{equation}\label{sd-tilde}
    \widetilde{s}(\mu)=\inf\big\{\gamma\in(0,1):\ \mu\ \textrm{is}\ \gamma\textrm{-spectral singular}\big\}.
\end{equation}
Obviously, ${s}(\mu)\le\widetilde{s}(\mu)$. The main theorem also holds for $\widetilde{s}(\mu)$.
\end{rmk}

\subsection{Spectral singularity, continuity and proof of Theorem \ref{thmmain}}
We first study  $\gamma$-spectral singularity of $\mu$. We are going
to show that under the assumption of Theorem \ref{thmmain} we have:
\begin{thm}\label{thmSingular}Assume $L(E)>a>0$ for $E\in S$. There
  exists $c=c(a)>0$ such that for any $\alpha,\theta,$
if
\begin{equation}\label{mSing}
    \gamma>\frac{1}{1+\frac{c}{\beta(\alpha)}},
\end{equation}
then $\mu_{S}$ is $\gamma$-spectral singular.
\end{thm}
Obviously, Theorem \ref{thmSingular} implies that if $\beta<+\infty$, then
\begin{equation}\label{upperDim}
    s(\mu_{S})\le \widetilde{s}(\mu_{S})\le \frac{1}{1+c/\beta}<1.
\end{equation}

 The analyticity of potential and positivity of Lyapunov exponent are
 only needed for spectral singularity. We now formulate a more precise
 version of the general spectral continuity
 result, Theorem \ref{00}.

For $S\subset \sigma(H)$ assume there are constants
 $\Lambda>0$ and $n_0\in \N$ such that for any 
$k\in\mathbb{Z},E\in S$ and $n\ge n_0$
\begin{equation}\label{Lambda}
\Big\|\left(\begin{array}{cc}
E- V(n+k)& -1 \\
1 & 0
\end{array}\right) \cdots \left(\begin{array}{cc}
E- V(k)& -1 \\
1 & 0
\end{array}\right)\Big\|\le e^{\Lambda n}
\end{equation}

Clearly, such $\Lambda$ always exists for bounded $V,$ with $n_0=1.$

As before we denote $\mu_{S}$  the spectral measure of $H$
restricted to a  Borel set $S\subset \sigma(H)$.
\begin{thm}\label{thmConti}Let $H$ be given by  (\ref{Schro-op-gen})
and  $V$ satisfies  (\ref{Lambda}) and is $\beta$-almost periodic with $\epsilon>0$. 
Then, for  a $C(\epsilon)=C_0(1+1/\epsilon)$, with
$\Lambda$ given by  (\ref{Lambda}) , if
\begin{equation}\label{mConti}
    \beta>C(\epsilon)\frac{\Lambda}{1-\gamma}
\end{equation}
then $\mu_{S}$ is $\gamma$-spectral continuous. Here $C_0$  is a universal constant.
Consequently, we have
\begin{equation}\label{lowerDim}
   \widetilde{s}(\mu_{S})\ge s(\mu_{S})\ge 1-C(\epsilon)\frac{\Lambda}{\beta}.
\end{equation}
\end{thm}

\noindent \textbf{Proof of Theorem \ref{thmmain}:}
Under the assumption of Theorem \ref{thmmain}, if $\beta<+\infty$,  Theorem \ref{thmSingular} provides the upper bound (\ref{upperDim}) for the spectral dimension.

We will now  get the lower bound using Theorem \ref{thmConti}. Let
$V_\theta(n):=V(\theta+n\alpha)$. By boundedness of $V$ and compactness of the spectrum,
 there is a constant $\Lambda_{V}<\infty$ 
such
that (\ref{Lambda}) holds uniformly for $E\in
\sigma(H_\theta),\theta\in\T$. In order to apply Theorem
\ref{thmConti}, it is enough to show that for any
$\beta<\beta(\alpha),$ $V(\theta+j\alpha)$ has $\beta$-repetitions for any $\theta\in\T,j\in\Z.$ 
Indeed, by (\ref{beta}), 
there is a subsequence $q_{n_k}$ such that $$\log q_{n_k+1}\big/q_{n_k}>\beta
$$
Since $V$ is analytic, for any $\theta,j$ and $1\le n\le q_{n_k}$
$$|V(\theta+j\alpha+n\alpha)-V(\theta+j\alpha+n\alpha\pm q_{n_k}\alpha)|\le C\|q_{n_k}\alpha\|\le C\frac{1}{q_{n_k+1}}\le Ce^{-\beta
  q_{n_k}}$$
Thus if $\beta(\alpha)=\infty$, $\widetilde{s}(\mu_{S})= s(\mu_{S})=1$. \qed


Property (\ref{Lambda}) naturally holds in a {\it sharp} way in the context of ergodic
potentials with uniquely ergodic underlying dynamics. Assume the potential $V=V_\theta$ is generated by some homeomorphism $T$ of a compact metric
space $\Omega$ and a function $f:\Omega\to \R$ by
\begin{equation}\label{Vtheta}
V_\theta(n)=f(T^n\theta), \ \theta\in\Omega,\ n\in\Z.
\end{equation}
Assume $(\Omega,T)$ is uniquely ergodic with an ergodic measure $\nu$. It is known that the spectral type
of $H_\theta$ is $\nu$-almost surely independent of $\theta$ (e.g \cite{cl}). In general,
however, the spectral type (locally) does depend on $\theta$ (\cite{js}).
If $f$ is continuous then, by uniform upper-semicontinuity (e.g. \cite{f97})
\begin{equation}\label{uniupper}
    \limsup_n\sup_\theta\frac{1}{n}\Big\|\left(\begin{array}{cc}
E- V_\theta(n)& -1 \\
1 & 0
\end{array}\right) \cdots \left(\begin{array}{cc}
E- V_\theta(1)& -1 \\
1 & 0
\end{array}\right)\Big\|\le L(E), \ \ \forall \ E
\end{equation}
This was recently extended in \cite{jmv2} to almost continuous
$f.$ 
Following \cite{jmv2}, we will say a function $f$ is almost continuous
if it is bounded and its set of discontinuities has a closure
of $\nu$ measure zero. By Corollary 3.2 in \cite{jmv2}, if $f$ is
bounded and almost continuous then (\ref{uniupper}) also holds for
every $E.$  
Moreover, if the Lyapunov exponent $L(E)$ is continuous on some
compact set $S$, then, by compactness and subadditivity,  the lim sup in (\ref{uniupper}) will be also
uniform in $E\in S$. Since by upper semicontinuity $L(E)$ is
continuous on the set where it is zero, as a consequence of Theorem
\ref{thmConti}, we obtain
\begin{cor}\label{zeroLyp}
Assume the function $f$ in (\ref{Vtheta}) is bounded and almost continuous and $L(E)=0$ on some Borel subset $S$ of $\sigma(H_\theta)$. If $V_{\theta}(n)$ is $\beta$-almost periodic  for some $\beta>0,\epsilon>0$, then 
$s(\mu^\theta_S)=1.$
\end{cor}
\proof For any $0<\gamma<1$, set $\Lambda'=\beta(1-\gamma)/2C$
where $C=C(\epsilon)$ is given in Theorem
\ref{thmConti}.\footnote{If $\beta=\infty$ take any finite $\beta$
  instead.} Since $L(E)=0$ on $S$, by the arguments above, there is $n_0=n_0(\Lambda')$ independent of $\theta$ and $E$ such that
$$\Big\|\left(\begin{array}{cc}
E- V_\theta(n)& -1 \\
1 & 0
\end{array}\right) \cdots \left(\begin{array}{cc}
E- V_\theta(1)& -1 \\
1 & 0
\end{array}\right)\Big\|\le e^{\Lambda' n},\ n\ge n_0,\ E\in S,\theta\in\Omega$$
Obviously, $\beta>C\Lambda'/(1-\gamma)$, so Theorem \ref{thmConti} is applicable and (\ref{lowerDim}) holds. Therefore,
$s(\mu_{S,\theta})\ge 1-\frac{C\Lambda'}{\beta}>\gamma.$ 
.\qed

Let $S_0=\{E:L(E)=0\}$ and $S_+=\{E:L(E)>0\}.$

As an immediate consequence we obtain
\begin{thm} \label{s0+}
If $V_\theta(n)$ is given by (\ref{Vtheta}) with uniquely ergodic
$(\Omega,T)$ and almost continuous $f,$ then, for every $\theta$ we
have
\begin{enumerate}
\item $s(\mu_{S_0})=1,$ as long as $V$ is $\beta$-almost periodic with $\beta>0.$
\item $s(\mu_{S_+})=1,$ as long as  $V$ is $\beta$-almost periodic with
  $\beta=\infty.$
\end{enumerate}
\end{thm}
\begin{rmk}
\begin{enumerate}\item
$\beta >0$ is not a necessary condition in general for
$s(\mu_{S_0})=1$, for $s(\mu^{ac})=1$ even if  $V$ is not
$\beta$-almost periodic for any $\beta$,  and the support
of the absolutely continuous spectrum is contained in (and may
coinside with) $S_0.$ It is a very interesting question to specify a quantitative
almost periodicity condition for $s(\mu^{sing}_{S_0})=1$, in particular, find
an arithmetic criterion  for
analytic one frequency potentials for $s(\mu_{S_{cr}})=1$ where
$S_{cr}\subset S_0$ is the set of critical energies in the sense of
Avila's global theory.
\item According to  Theorem \ref{thmmain}, $\beta =\infty$ is also
  necessary if $f$ is analytic and $T$ is an irrational rotation of
  the circle ($\beta$ will depend on $T$). In case $f$ has lower
  regularity, it is an interesting  question to determine optimal
  condition on $\beta$.
\end{enumerate}\end{rmk}
\subsection{Relation with other dimensions; Corollaries for the AMO,
  Sturmian potentials, and Transport exponents.}\label{secApp}
If we replace the $\liminf$ by $\limsup$ in the definition of upper
spectral dimension, we will define correspondingly the lower spectral
dimension which will coincide with the Haurdorff dimension $dim_H(\mu)$ of a measure $\mu$.

Also one can consider the packing dimension of $\mu$, denoted by
$dim_P(\mu)$. The packing dimension can be defined in a similar way as
in (\ref{sd}) through the $\gamma$-dimensional lower derivative
$\underline{D}^{\gamma}\mu(E)$. It can be easily shown that
$\underline{D}^{\gamma}\mu(E)\le
\liminf_{\varepsilon\downarrow0}\varepsilon^{1-\gamma}|M(E+ i\varepsilon)|$.
Thus the relation between packing dimension and upper spectral dimension is
$dim_P(\mu)\ge \widetilde{s}(\mu)$.\footnote{In contrast with the Hausdorff
  dimension, the relation for the packing dimension only goes in one
  direction, in general, unlike what is claimed in \cite{co1}.} Therefore, the lower bound we get
in Theorem \ref{thmConti} also holds for the packing dimension.

Lower bounds on  spectral dimension also have immediate applications to
 the lower bounds on quantum dynamics. Denote by $\delta_j$ be the vector $\delta_j(n)=\chi_{j}(n)$. For $p > 0$, define
\begin{equation}\label{mom}
    \langle|X|_{\delta_0}^p\rangle(T)=\frac{2}{T}\int_0^\infty e^{-2t/T}\sum_n|n|^p|\langle e^{- itH}\delta_0,\delta_n\rangle|^2
\end{equation}
The growth rate of $\langle|X|_{\delta_0}^p\rangle(T)$ characterizes how fast does $e^{- itH}\delta_0$ spread out. In order to get the power law bounds for $\langle|X|_{\delta_0}^p\rangle(T)$, it is natural to define the following upper $\beta_{\delta_0}^+(p)$ and lower $\beta_{\delta_0}^-(p)$ dynamical exponents  as
\begin{equation}\label{dyn-exp}
    \beta_{\delta_0}^+(p)=\limsup_{T\to\infty}\frac{\log\langle|X|_{\delta_0}^p\rangle(T)}{p\log T}, \ \
    \beta_{\delta_0}^-(p)=\liminf_{T\to\infty}\frac{\log\langle|X|_{\delta_0}^p\rangle(T)}{p\log T}
\end{equation}
The dynamics
is called ballistic if $\beta_{\delta_0}^-(p)=1$ for all $p>0$, and
quasiballistic if $\beta_{\delta_0}^+(p)=1$ for all $p>0$. We will
also say that the dynamics is quasilocalized if $\beta_{\delta_0}^-(p)=0$ for all $p>0$.


In \cite{gsb}, it is shown that the upper and lower transport exponents
of a discrete Schr\"odinger operator (\ref{Schro-op-gen}) can be
bounded from below by the packing and Hausdorff dimension of its
spectral measure respectively. Therefore, by \cite{gsb} we have $\beta_{\delta_0}^+(p)\ge s(\mu), \ \forall p$.  As a direct consequence of Theorem \ref{thmConti} 
we have
\begin{cor}
If $V(n)$ is bounded  and $\infty$-almost periodic, 
the upper dynamical exponent $\beta_{\delta_0}^+(p)$ of the operator (\ref{Schro-op}) is one for any $p>0$, and the associated dynamics is quasiballistic.

\end{cor}

This has nice immediate consequences. In
particular, consider the almost Mathieu operator:
\begin{equation}\label{AMO}
    (H_{\lambda,\theta,\alpha} u)_n= u_{n+1}+u_{n-1} + 2\lambda\cos2 \pi
(\theta+n\alpha)u_n,\ \lambda>0.
\end{equation}
As a consequence of the formula for the Lyapunov exponent and Theorem \ref{thmmain}, one has:
\begin{cor}\label{corAMO}The
almost Mathieu operator (\ref{AMO}) is quasiballistic\footnote{And
  has spectral dimension one and packing dimension one of the spectral
  measure} for any  (and all) $\theta\in\T$
\begin{enumerate}
  \item For $\lambda<1,$ for all $\alpha$
  \item For $\lambda=1,$ as long as $\beta(\alpha)>0$
  \item For $\lambda>1,$ as long as $\beta(\alpha)=\infty$.
\end{enumerate}
\end{cor}

Statement 1 is a corollary of absolute continuity \cite{last93,a-ac} and
is listed here for completeness only.
 Statements 2,3 are direct corollaries of Theorem \ref{s0+}.

For $\lambda>1,$ Hausdorff
dimension of the spectral measure of the almost Mathieu operator  is equal to zero \cite{jl2} and
$\beta^-(p)=0$ for all $p>0$ \cite{dt}.  Thus almost Mathieu operators
with $\lambda>1$ and $\beta(\alpha)=\infty$ provide a family of
{\it explicit} examples of operators that are simultaneously quasilocalized
and quasiballistic and whose spectral measures satisfy
$$0=\dim_H(\mu)<\dim_P(\mu)=1.$$ The same holds of course for $\cos$
replaced with any almost continuous $f$ as long as Lyapunov exponent
is positive everywhere on the spectrum, in particular for $f=\lambda
g$  where $g$ is either bi-Lipshitz (as in \cite{jk}) or analytic,
and $\lambda>\lambda(g)$ is sufficiently large.

Let ${\rm d}N$ be the density states measure of the almsot Mathieu opeartor and $\Sigma$ be
the spectrum. It is well known that in the critical case, $\lambda=1$,
$\Sigma$ has Lebesgue measure zero (\cite{ak,last94}). It is then
interesting to consider the fractal dimension of the spectrum (as a
set). Since ${\rm d}N=\E({\rm d}\mu_\theta)$ and ${\rm
  supp_{top}}({\rm d}N)= \Sigma$ , by the discussion above we have
\begin{cor}\label{corAMO1}
For the critical almost Mathieu operator, $\lambda=1,$ and $\beta(\alpha)>0$ we have
$dim_P(dN)=dim_P(\Sigma)=1.$
\end{cor}

Last and Shamis proved in \cite{ls} (see also \cite{mirathesis}) that for a dense $G_\delta$ set of
$\alpha$ (which therefore has a generic intersection with the set
$\{\alpha: \beta(\alpha)>0\}$), the Hausdorff dimension of the
spectrum is equal to zero. Thus the spectrum of the critical
almost Mathieu operator with a topologically
generic frequency is an example of a {\it set} such that
$$0=\dim_H(\Sigma)<\dim_P(\Sigma)=1.$$ Moreover, Last \cite{last94} showed that if
$q_{n+1}>Cq_n^3$ for all $n,$ (which is the set containing $\{\alpha:
\beta(\alpha)>0\}$) then $dim_H(\Sigma)\leq 1/2$. Thus critical almost
Mathieu operators with $\beta(\alpha)>0 $ and any $\theta$ provide an
explicit family of operators that all have
spectra satisfying $\dim_H(\Sigma)\leq 1/2<\dim_P(\Sigma)=\dim_B(\Sigma)=1.$ \\

We note that the question of fractal dimension of the critical almost
Mathieu operator attracted a lot of attention in Physics literature,
with many numerical and heuristic results. In particular, Wilkinson-Austin \cite{wa}
conjectured that $\dim_B(\Sigma)<1/2$ for {\it all} $\alpha$ and there were
many results rigorously or numerically confirming this for
certain $\alpha$. Our corollary \ref{corAMO1} provides an explicit
example disproving this conjecture.\\

Another  well known family are  Sturmian Hamiltonians given by
\begin{equation}\label{defSturm}
    (H u)_n= u_{n+1}+u_{n-1} + \lambda \chi_{[1-\alpha,1)}(n \alpha + \theta \ \textrm{mod}\ 1)u_n,
\end{equation}
where $ \lambda > 0, \ \alpha = \R\backslash\Q$. If $\alpha =
\frac{\sqrt{5}-1}{2}$, it is called the Fibonacci Hamiltonian. The
spectral properties of the Fibonacci Hamiltonian have been thoroughly
studied in a series of papers in the past three decades, see \cite{d07,degt} for more references. Recently, Damanik, Gorodetski, Yessen proved in \cite{dgy} that for every $\lambda > 0$, the density of states measure $\textrm{d}N_\lambda$ is exact-dimensional (the Hausdorff and upper box counting dimension are the same) and ${\rm dim}_H(\textrm{d}N_\lambda)<{\rm dim}_H(\Sigma_\lambda)$.

Our results show that the exact dimensionality properties of Sturmian
Hamiltonians strongly rely on the arithmetic properties of
$\alpha$. It was shown in \cite{bist} that if $\alpha$ is irrational,
the Lyapunov exponent of Sturmian operator restricted to the spectrum is
zero. Also the spectrum of Sturmian Hamiltonian
$\Sigma_{\lambda,\alpha}$ is always a Cantor set with Lebesgue measure
zero. Moreover, for Sturmian potentials results similar
to thoses for the critical Mathieu operator in Corollary \ref{corAMO1}
also hold.
Let $\mu_{\theta}$ be the spectral measure of Sturmian operator
(\ref{defSturm}) and let $\textrm{d}N_{\lambda,\alpha}$ be the density
states of measure and $\Sigma_{\lambda,\alpha}$ be the spectrum. We
say that phase $\theta$ is $\alpha$-Diophantine if there exist
$\gamma<\infty,\tau>1$ such that $\|\theta+m\alpha\|_{\R/\Z}\geq
\frac{\gamma}{(|m|+1)^\tau}$ for all $m\in\mathbb{Z}.$ Clearly, this
is a full measure condition. We have
\begin{thm}\label{thmSturm}
For Sturmian operator $H_{\theta,\lambda,\alpha}$ with $\beta(\alpha)>0$ and $\lambda>0,$ if $\theta$ is $\alpha$-Diophantine, the spectral dimension of $\mu_{\theta}$ is one.

 As a consequence, if $\beta(\alpha)>0$ and $\lambda>0$, then the packing dimension of  $\textrm{d}N_{\lambda,\alpha}$ and $\Sigma_{\lambda,\alpha}$  are both equal to one.
\end{thm}

Previously, Liu, Qu and Wen \cite{lw04,lqw14} studied the
Hausdorff and upper box counting dimension of
$\Sigma_{\lambda,\alpha}$ of Sturmian operators. For large couplings, they gave
a criterion on $\alpha\in(0,1)$ for the Hausdorff dimension of the spectrum to
be equal to one. Combining Theorem \ref{thmSturm} with their results, we have

\begin{cor}\label{corSturm}
Let $\Sigma_{\lambda,\alpha}$  be the spectrum of the Sturmian
Hamiltonian with $\lambda>20$. There are explicit $\alpha$ such that
for a.e. $\theta$,
\begin{equation}
    {\rm dim}_H(\mu^\theta_{\lambda,\alpha})<s(\mu^\theta_{\lambda,\alpha})={\rm dim}_P(\mu^\theta_{\lambda,\alpha})=1
\end{equation}
\begin{equation}
    {\rm dim}_H(\textrm{d}N_{\lambda,\alpha})<s(\textrm{d}N_{\lambda,\alpha})={\rm dim}_P(\textrm{d}N_{\lambda,\alpha})=1
\end{equation}
\end{cor}

The proof for the Sturmian case is given in Section \ref{secSturm}.\\

The rest of this paper is organized in the following way. After giving
the preliminaries in Section  \ref{prel} we proceed to the proof of
the general continuity statement in Section \ref{secConti}. First we
quickly reduce Theorem \ref{thmConti} to Lemma \ref{lemPowerlaw}  where we also
specify the constant $C_0$ appearing in Theorem \ref{thmConti}. We note that
we do not aim to optimize the constants here and many of our arguments have
room for corresponding improvement. Lemma \ref{lemPowerlaw} is further reduced
to the estimate on the traces of the transfer matrices over eventual
almost periods, Theorem \ref{thmtrace}, through its corollaries,
Lemmas \ref{lem4} and \ref{lem41}. Theorem \ref{thmtrace} is the  key
element and the most technical part of the proof.  It is of interest
in its own right as can be viewed as the quantitative version of the
fact that period length transfer matrices of periodic operators are
elliptic: it provides quantitative bounds on the traces of transfer
matrices over almost periods based on
quantitative almost periodicity, for spectrally a.e. energy. In section
\ref{sec2.2} we separate this statement into hyperbolic and almost
parabolic parts, correspondingly Lemmas \ref{lem1} and \ref{lem2}.  In
section \ref{sec2.3} we use the extended Schnol's  Theorem to study
the hyperbolic case and in section \ref{sec2.4} we combine estimates on
level sets of the polynomials, power-law subordinacy bounds, and an
elementary but very useful
algebraic representation of matrix powers (Lemma \ref{lemAkqFormula})  to study the almost parabolic
case. Lemmas \ref{lem4} and \ref{lem41} are proved in Section \ref{sectionlem4},
completing the continuity part. In Section \ref{secSing} we  focus on the
analytic quasiperiodic potentials and prove Theorem
\ref{thmSingular}. The proof is  based on a lemma about density of
localized blocks (Lemma \ref{lemexp}). Finally, we discuss
Sturmian potentials in Section \ref{secSturm}, proving Theorem
\ref{thmSturm} and then providing explicit examples  for Corollary
\ref{corSturm}.

\subsection{Preliminaries}\label{prel}
\subsubsection{m-function and subordinacy theory}
In this part, we will briefly introduce the power-law extension of the
Gilbert-Pearson subordinacy theory
\cite{g,gp}, developed in \cite{jl1}. We will also list the necessary
related facts on the Weyl-Titchmarsh $m$-function. More details can be found, e.g., in \cite{cl}.\\

Let $H$ be as in (\ref{Schro-op}) and
$z=E+ i\varepsilon\in\C.$ Consider equation
\begin{equation}\label{eig-eq}
    Hu=zu.
\end{equation}
with the family of normalized phase boundary conditions:
\begin{equation}\label{bdry}
    u^{\varphi}_0\cos\varphi+u^{\varphi}_1\sin\varphi=0,\ -\pi/2<\varphi<\pi/2,\  |u^{\varphi}_0|^2+|u^{\varphi}_1|^2=1.
\end{equation}
Let $\Z^+=\{1,2,3\cdots\}$ and $Z^-=\{\cdots,-2,-1,0\}$. Denote by $u^{\varphi}=\{u^{\varphi}_j\}_{j\ge0}$ the right half line solution on $\Z^+$ of (\ref{eig-eq}) with  boundary condition (\ref{bdry}) and by $u^{\varphi,-}=\{u^{\varphi,-}_j\}_{j\le0}$ the left half line solution on $\Z^-$ of the same equation. Also denote by $v^{\varphi}$ and $v^{\varphi,-}$  the right and left half line solutions of (\ref{eig-eq}) with the orthogonal boundary
conditions to $u^{\varphi}$ and $u^{\varphi,-}$, i.e., $v^{\varphi}=u^{\varphi+\pi/2}$,$v^{\varphi,-}=u^{\varphi+\pi/2,-}$. For any function $u:\Z^+\to\C$ we denote by $\|u\|_l$ the norm of $u$ over a lattice
interval of length $l$; that is
\begin{equation}\label{ul}
\|u\|_l=\Big[\sum_{n=1}^{[l]}|u(n)|^2+(l-[l])|u([l]+1)|^2\Big]^{1/2}
\end{equation}
Similarly, for $u:\Z^-\to\C$, we define
\begin{equation}\label{ul-}
\|u\|_l=\Big[\sum_{n=1}^{[l]-1}|u(-n)|^2+(l-[l])|u(-[l])|^2\Big]^{1/2}
\end{equation}
Now given any $\varepsilon>0$, we define lengths $l=l(\varphi,\varepsilon,E)$, by requiring the equality
\begin{equation}\label{l}
    \|u^{\varphi}\|_{l(\varphi,\varepsilon)}\|v^{\varphi}\|_{l(\varphi,\varepsilon)}=\frac{1}{2\varepsilon}
\end{equation}
We also define $l^-(\varphi)$ by $u^{\varphi,-},v^{\varphi,-}$ through the same equation.
Direct computation shows that
\begin{equation}\label{uv0}
   \|u^{\varphi}\|_{l}\cdot \|v^{\varphi}\|_{l}\ge \frac{1}{2}([l]-1)
\end{equation}

Denote by $m_{\varphi}(z):\C^+\mapsto\C^+$ and
$m_{\varphi}^-(z):\C^+\mapsto\C^+$ the right and left Weyl-Tichmarsh
m-functions associated with the boundary condition (\ref{bdry}). Let
$m=m_0$ and $m^-=m^-_0$ be  the half line m-functions corresponding to the Dirichlet boundary conditions.
The following key  inequality \cite{jl1} relates $m_{\varphi}(E+ i\varepsilon)$ to the solutions $u^{\varphi}$ and $v^{\varphi}$ given by (\ref{eig-eq}),(\ref{bdry}).
\begin{lemma}[J-L inequality, Theorem 1.1 in \cite{jl1}] \label{lemJL}
For $E\in\R$ and $\varepsilon>0$, the following inequality holds for any $\varphi\in(-\frac{\pi}{2},\frac{\pi}{2}]$:
\begin{equation}\label{jl}
    \frac{5-\sqrt{24}}{|m_{\varphi}(E+i\varepsilon)|}
<\frac{\|u^{\varphi}\|_{l(\varphi,\varepsilon)}}{\|v^{\varphi}\|_{l(\varphi,\varepsilon)}}
<\frac{5+\sqrt{24}}{|m_{\varphi}(E+i\varepsilon)|}
\end{equation}

\end{lemma}

We need to study the whole-line m-function which is given by the Borel transform of the spectral measure $\mu$ of operator $H$ (see e.g., \cite{cl}).
The following relation between whole line m-function $M$ and half line m-function $m_\varphi$ was first shown
in \cite{dkl} as a corollary of the maximal modulus principle. One can
also find a different proof based on a direct computation in the
hyperbolic plane in \cite{aj}.

\begin{prop}[Corollary 21 in \cite{dkl}] \label{DKL}
Fix $E\in\R$ and $\varepsilon>0$,
\begin{equation}\label{aj2}
    |M(E+ i\varepsilon)|\le\sup_{\varphi}|m_{\varphi}(E+ i\varepsilon)|
\end{equation}
\end{prop}
This proposition implies that in order to obtain an upper bound for
the whole line m-function, namely, the continuity of whole line
spectrum, it is enough to obtain a uniform upper bound of the half line m-function  for any boundary condition.

On the other hand, consider a unitary operator $U: l^2(\Z)\to
l^2(\Z)$, defined by $(U\psi)_n=\psi_{-n+1},\ n\in\Z$.  For any
operator $H$ on $l^2(\Z)$, we define an operator $\widetilde{H}$ on
$l^2(\Z)$ by $\widetilde{H}=UHU^{-1}$. Denote by
$\widetilde{m},\widetilde{m}_\varphi,\widetilde{u}^{\varphi}$ and
$\widetilde{l}(\varphi)$, correspondingly,
${m},{m}_\varphi,{u}^{\varphi}$ and ${l}(\varphi)$ of the operator
$\widetilde{H}$. We will need the following  well known facts (see e.g. Section 3, \cite{jl2}). For any $\varphi\in(-\pi/2,\pi/2]$ we have
\begin{equation}\label{Mm+m-}
    M(z)=\frac{m_\varphi(z)\widetilde{m}_{\pi/2-\varphi}-1}{m_\varphi(z)+\widetilde{m}_{\pi/2-\varphi}}
\end{equation}
and
\begin{equation}\label{L+-}
    \widetilde{l}(\pi/2-\varphi)={l}^-(\varphi),\  \ \|u\|_l=\|Uu\|_l
\end{equation}

Similar to Lemma 5 in \cite{jl2}, a direct consequence of relation (\ref{Mm+m-}) is the following result.
\begin{lemma}\label{lemMmm}
For any $0<\gamma<1$, suppose that there exists a $\varphi\in(-\pi/2,\pi/2]$ such that for $\mu$-a.e. $E$ in some Borel set $S$, we have that $\liminf_{\varepsilon\to0} \varepsilon^{1-\gamma}|m_\varphi(E+ i\varepsilon)|=\infty$ and $\liminf_{\varepsilon\to0} \varepsilon^{1-\gamma}|\widetilde{m}_{\pi/2-\varphi}(E+ i\varepsilon)|=\infty$.
Then for $\mu$-a.e. $E$ in $S$, $\liminf_{\varepsilon\to0} \varepsilon^{1-\gamma}|M(E+ i\varepsilon)|=\infty$, namely,  the restriction $\mu(S\bigcap\cdot)$ is $\gamma$-spectral singular.
\end{lemma}

\subsubsection{Transfer matrices and Lyapunov exponents}\label{pre2}
Although Theorem \ref{thmConti} does not involve any further
conditions on the potential, it will be convenient in what follows to
use the dynamical notations. Let $\Omega=\R^\Z$ and
$T:\Omega\mapsto\Omega$ is given by $(T\theta)(n)=\theta(n+1)$.
Let  $f(\theta):=\theta(0).$ Then any potential $V$ can be written  in the way  (\ref{Vtheta}), $V_\theta(n):=\theta(n)=f(T^n\theta)$.
Thus for a fixed $\{V_n\}_{n\in\Z}=\theta\in\Omega$, we will rewrite  the
potential $V$ as $V_\theta(n)=f(T^n\theta)$ as in (\ref{Vtheta}). For
our general theorem we do not introduce any topology, etc; this is
being done purely for the notational convenience.   Denote the $n$-step transfer-matrix by $A_n(\theta,E)$:
\begin{equation}\label{ntransfer}
    A_n(\theta,E)=
A\big(T^n\theta,E\big)A\big(T^{n-1}\theta,E\big)\cdots A\big(T\theta,E\big),\ \ n>0
\end{equation}
and
$$A_0=Id,\ \ \ A_{n}(\theta,E)=A_{-n}^{-1}(T^n\theta,E),\ \ n<0,$$
where
\begin{equation}\label{Schro-cocycle}
    A(\theta,E)=\left(\begin{array}{cc}
E- f(\theta)& -1 \\
1 & 0
\end{array}\right).
\end{equation}

The connection to Schr\"{o}dinger operators is clear since a solution of $Hu=Eu$ can be reformulated as
\begin{equation}\label{ref-sln}
    A_{n}(\theta,E)\left(
  \begin{array}{c}
    u_{1} \\
    u_{0} \\
  \end{array}
\right)=\left(
  \begin{array}{c}
    u_{n+1} \\
    u_{n} \\
  \end{array}
\right), \ n\in\Z .
\end{equation}
In other words, the spectral properties of Schr\"{o}dinger operators $H$
are closely related to the dynamics of the family of skew product
$(T,A(\theta,E))$ over $\Omega\times\R^2$. We will often suppress
either $\theta$ or $E$ or both from the notations if corresponding
parameters are fixed through the argument.

If $V$ is
actually dynamically defined by  (\ref{Vtheta}) with a certain
underlying ergodic base dynamics $(\Omega,T,\nu)$ then, by the general properties of subadditive ergodic cocycles, we can define the Lyapunov exponent
\begin{equation}\label{defLyp}
    L(E,T)=\lim_{n\to +\infty}\frac{1}{n}\int_\Omega\log\|A_n(\theta,E)\|{\rm d}\nu=
    \inf_{n>0}\frac{1}{n}\int_\Omega\log\|A_n(\theta,E)\|{\rm d}\nu
\end{equation}



\section{Spectral Continuity }\label{secConti}
\subsection{Proof of Theorem \ref{thmConti}}
Throughout this section 
we assume (\ref{Lambda}) is satisfied  uniformly for  $E\in
S,\theta\in\Omega$ (see Section 1.4.2).  
Assume $V$ is
$\beta$-almost periodic for some $\epsilon>0$. 
The proof of Theorem \ref{thmConti} is based on the following estimates on the growth of the $l$-norm of the half line solutions. Let $u^{\varphi},v^{\varphi}$ be given as in (\ref{eig-eq})-(\ref{ul}).
\begin{lemma}\label{lemPowerlaw}
For $0<\gamma<1,$ assume $\beta>100(1+1/\epsilon)\Lambda/(1-\gamma)$.
For $\mu$-a.e. $E$, 
there is a sequence of positive numbers $\eta_k\to 0$ so that
for any $\varphi$ 
\begin{equation}\label{uvLk}
  1/16\big(L_k\big)^{\gamma} \le \|v^{\varphi}\|^2_{L_k}\le  \big(L_k\big)^{2-\gamma}
\end{equation}
where $L_k=\ell_k(\varphi,\eta_k,E)$ is given as in (\ref{l}).
\end{lemma}
\noindent
\textbf{Proof of Theorem \ref{thmConti}:}
Fix $0<\gamma<1$. Set $C(\epsilon):=
100(1+1/\epsilon)$.  Lemma \ref{lemPowerlaw} can be applied to any $\beta>C(\epsilon)\frac{\Lambda}{1-\gamma}$.
According to (\ref{uvLk}) and the  J-L inequality (\ref{jl}), for $\mu$-a.e. $E$ and any $\varphi$
\begin{eqnarray*}
  \eta_k^{1-\gamma}|m_\varphi(E+ i\eta_k)| &\le&
\frac{1}{\Big(2 \|u^{\varphi}\|_{L_k}\|v^{\varphi}\|_{L_k}\Big)^{1-\gamma}}\cdot (5+\sqrt24)\frac{\|v^{\varphi}\|_{L_k}}{\|u^{\varphi}\|_{L_k}}
 \\
&\le&
C_{\gamma}\cdot \frac{\big( L_k^{(2-\gamma)/2}\big)^{\gamma}}{\big(
  1/4L_k^{\gamma/2}\big)^{2-\gamma}} = 
C <\infty
\end{eqnarray*}

Since $\eta_k$ is independent of $\varphi$, for a fixed $E$ and
$\eta_k$, we can take the supremum w.r.t. $\varphi$. By (\ref{aj2}) in Proposition \ref{DKL}, we have for $\mu$-a.e. $E$,
$$
    \eta_k^{1-\gamma}|M(E+ i\eta_k)|<C
$$
i.e.,
$$
    \liminf_{\varepsilon\downarrow0}\varepsilon^{1-\gamma}|M(E+ i\varepsilon)|<\infty, \ \ \mu\textrm{-a.e.}\ E,
$$
which proves the $\gamma$-spectral continuity of  Theorem
\ref{thmConti}. 
The lower bound (\ref{lowerDim}) comes from the definition of spectral dimensionality.  \qed

The proof of Lemma \ref{lemPowerlaw} follows from the following
estimates on the trace of the transfer matrix. 
Let $q_k$ be the sequence given in (\ref{beta-exp-peri}).
\begin{thm}\label{thmtrace}
If
\begin{equation}\label{betathmtrace}
    \beta>(37+11/\epsilon)\Lambda,
\end{equation}
then for $\mu$ a.e. $E$, there is $K(E)$ such that
\begin{equation}\label{trace}
    |{\rm Trace}A_{q_k}(E)|<2-e^{-10\Lambda q_k}, \ \ k\ge K(E).
\end{equation}

\end{thm}
This theorem is the key estimate of the spectral continuity. It can be
viewed as a quantitative version of the classical fact that
period-length  transfer
matrices of periodic operators are elliptic on the spectrum. Indeed,
we prove that $\beta$-almost periodicity implies quantitative bounds
on ellipticity.
The proof
will be given in the following two subsections. The direct consequence
of Theorem \ref{thmtrace} are the following estimates on the norm of the
transfer matrices. They show that if the trace of the transfer matrix
over an almost period
is strictly less than $2$, then the there is a sub-linearly
bounded subsequence. We will use this result to prove Lemma
\ref{lemPowerlaw} first. The proof of Lemma \ref{lem4} will be left to
Section \ref{sectionlem4}. Let $K(E)$ be given by Theorem \ref{thmtrace}.
\begin{lemma}\label{lem4}
For any $\xi>0$ set $N_k=[e^{\xi q_k}]$  and suppose that, in addition
to the conditions of Theorem \ref{thmtrace}
\begin{equation}\label{xi1}
    \beta>15\Lambda+(2+1/\epsilon)\xi.
\end{equation}
Then for $\mu\textrm{-a.e.}\ E$, the following estimate holds:
\begin{equation}\label{grow1}
    \sum_{n=1}^{N_k\cdot q_k}\|A_n(E)\|^2\le  e^{(\xi+15\Lambda)q_k}, \ \ k\ge K(E)
\end{equation}
\\
\end{lemma}
Additionally,
\begin{lemma}\label{lem41}
For $0<\gamma<1$, assume that in addition to the conditions of Lemma \ref{lem4}
\begin{equation}\label{xi2}
    \xi>\frac{16\Lambda}{1-\gamma}.
\end{equation}
Then
\begin{equation}\label{grow2}
 \sum_{n=1}^{N_k\cdot q_k}\|A_n(E)\|^2\le  (N_k\cdot q_k)^{2-\gamma}, \ \ k\ge K(E).
\end{equation}

\end{lemma}

\noindent \textbf{Proof of Lemma \ref{lemPowerlaw}:}
It is enough to prove r.h.s. of (\ref{uvLk}) since then the l.h.s. of (\ref{uvLk}) follows from
$\|u^{\varphi}\|_{L_k}\|v^{\varphi}\|_{L_k}\ge 1/4L_k$.


For any $0<\gamma<1$, set
$\beta_0=100(1+1/\epsilon)\frac{\Lambda}{1-\gamma}$,
$\xi=\frac{17\Lambda}{1-\gamma}$. Then (\ref{betathmtrace}),(\ref{xi1})
and (\ref{xi2}) are satisfied for all $\beta>\beta_0$. Therefore,
(\ref{grow2}) holds with above choice of parameters. Let $l_k=[e^{\xi
  q_k}]\cdot q_k.$  Rewrite (\ref{grow2}) as $
    \sum_{n=1}^{l_k}\|A_n(E)\|^2<l_k^{2-\gamma}$.
Thus for any $\varphi,$ $\|v^{\varphi}\|^2_{l_k}\le 4 l_k^{2-\gamma}$. By (\ref{uv0}),  we have
\begin{equation}\label{uv}
    \frac{1}{4}\cdot l_k\le\|u^{\varphi}\|_{l_k}\|v^{\varphi}\|_{l_k}\le 4\cdot l_k^{2-\gamma}
\end{equation}
Set
\begin{equation}\label{epsilon}
 \varepsilon_k(\varphi):= \frac{1}{2 \|u^{\varphi}\|_{l_k}\|v^{\varphi}\|_{l_k}}
\end{equation}
Then,
\begin{equation}\label{eta}
    \eta_k=\inf_{\varphi}\ \varepsilon_k(\varphi)\ge \frac{1}{8\cdot l_k^{2-\gamma}}>0
\end{equation}
is well defined. Set $L_k(\varphi):=\ell(\varphi,\eta_k,E),$ the length scale satisfies
\begin{equation}\label{Lk}
    \eta_k= \frac{1}{2 \|u^{\varphi}\|_{L_k(\varphi)}\cdot \|v^{\varphi}\|_{L_k(\varphi)}}
\end{equation}

By (\ref{eta}),
$$ L_k(\varphi)\le 4 \|u^{\varphi}\|_{L_k}\|v^{\varphi}\|_{L_k}=\frac{2}{\eta_k}\le 16\cdot l_k^{2-\gamma} $$
Since $\eta_k\le\varepsilon_k(\varphi)$ and $\|u^{\varphi}\|_{l}\|v^{\varphi}\|_{l}$ is monotone increasing in $l$, we obtain for any $\varphi$
\begin{equation}\label{Lk0}
   l_k\le L_k(\varphi)\le 16\cdot l_k^{2-\gamma}
\end{equation}

By the definition of $l_k$, for large $k,$
\begin{equation}\label{Lk1}
    e^{(\xi-\frac{\Lambda}{200(1-\gamma)}q_k}\cdot q_k\le L_k(\varphi)
\le  e^{((2-\gamma)\xi+\Lambda/200)q_k}\cdot q_k
\end{equation}

Write $L_k(\varphi)=[L_k(\varphi)]+\widetilde{L}_k(\varphi)$ and
\begin{equation}\label{Lk2}
    [L_k(\varphi)]=(N_k(\varphi)-1)\cdot q_k+r_k(\varphi), \ \ N_k(\varphi)\in \N, \ \ 0\le r_k(\varphi)< q_k
\end{equation}
where $[L_k(\varphi)]$ is the integer part of $L_k(\varphi)$. Define
\begin{equation}\label{bk}
    \xi_k(\varphi)=\frac{\log N_k(\varphi)}{q_k}
\end{equation}

We have
\begin{equation}\label{Lk3}
    [L_k(\varphi)]=(e^{\xi_k(\varphi)\cdot q_k}-1)\cdot q_k+r_k(\varphi), \ \ e^{\xi_k(\varphi)\cdot q_k}\in \N, \ \ 0\le r_k(\varphi)< q_k
\end{equation}
For large $q_k$, it is easy to check
$$e^{(\xi_k(\varphi)-\Lambda/200)\cdot q_k}\cdot q_k\le
(e^{\xi_k(\varphi)\cdot q_k}-1)\cdot q_k\le L_k(\varphi)\le e^{\xi_k(\varphi)\cdot q_k}\cdot q_k$$

Using (\ref{Lk1}), we have for any $\varphi$
\begin{equation}\label{bk1}
    \xi-\Lambda/200\le \xi_k(\varphi)\le (2-\gamma)\xi+\Lambda/100\le 2\xi+\Lambda/100
\end{equation}
Together with the choice of $\beta$ and $\xi$, we have
$$ \beta>\beta_0>15\Lambda+(2+1/\epsilon)\xi_k(\varphi)$$
Now we can again apply Lemma \ref{lem4} with parameters $\beta$, $\xi_k(\varphi)$ and the length scale
 $N_k(\varphi)=e^{\xi_k(\varphi)\cdot q_k}$, to get
\begin{equation}\label{Nk}
  \sum_{n=1}^{N_k(\varphi)\cdot q_k}\|A_n(E)\|^2\le e^{\big(\xi_k(\varphi)+15\Lambda\big)q_k}
\end{equation}

Notice $L_k(\varphi)\ge 
 e^{(\xi_k(\varphi)-\Lambda/200)\cdot q_k}$ implies that
\begin{eqnarray*}
   \frac{1}{(L_k)^{2-\gamma}}\sum_{n=1}^{N_k(\varphi)\cdot q_k}\|A_n(E)\|^2 &\le&
e^{\big(-(1-\gamma)\xi_k+16\Lambda\big)q_k}
\end{eqnarray*}

By the l.h.s. of (\ref{bk1}), we obtain
$$(1-\gamma)(\xi_k(\varphi)+\Lambda/200)>(1-\gamma)\xi=17\Lambda$$
which implies
$$
    (1-\gamma)\xi_k(\varphi)>17\Lambda-(1-\gamma)\Lambda/200>16.5\Lambda,
$$
and
\begin{equation}\label{Lk4}
    \frac{1}{(L_k)^{2-\gamma}}\sum_{n=1}^{N_k\cdot q_k}\|A_n(E)\|^2\le e^{-\Lambda q_k/2}\le 1
\end{equation}

Finally, by Lemma \ref{lem41} we have
$$
   \|v^{\varphi}\|_{L_k}^2    \le  \sum_{n=1}^{[L_k]+1}|v^{\varphi}_n|^2
   \le \sum_{n=1}^{N_k(\varphi)\cdot q_k}\Big(|v^{\varphi}_{n}|^2+|v^{\varphi}_{n+1}|^2\Big)
   \le  \sum_{n=1}^{N_k(\varphi)\cdot q_k}\|A_n(E)\|^2 \le (L_k)^{2-\gamma}
$$
\qed


\subsection{Proof of Theorem \ref{thmtrace}}\label{sec2.2}
The proof of Theorem \ref{thmtrace} will be divided into two cases. We
will first exclude the energies where the trace is much greater than 2
infinitely many times using extended Schnol's Theorem (Lemma
\ref{schnol}). Then we will estimate the measure of energies where the
trace is close to 2 through subordinacy theory. The conclusion
consists of the following two lemmas. Again let $q_k$ be the sequence given by (\ref{beta-exp-peri}) with certain $\beta,\epsilon>0$.
\begin{lemma}\label{lem1}
For any $\tau>0$, if
\begin{equation}\label{tau0}
    \beta>(3+1/\epsilon)\tau+(7+1/\epsilon)\Lambda,
\end{equation}
then for spectrally a.e. $E$, there is $K_1(E)$ such that,
\begin{equation}\label{K1}
    |{\rm Trace}\, A_{q_k}\big(E\big)|<2+e^{- \tau q_k}, \ \ \forall k\ge K_1(E)
\end{equation}

 \end{lemma}

\begin{lemma}\label{lem2}
If
\begin{equation}\label{lem2c2}
    \beta>(25+1/\epsilon)\Lambda,
\end{equation}
then for spectrally a.e. $E$, there is $K_2(E)$ such that
\begin{equation}\label{K2}
    |{\rm Trace}\,A_{q_k}(E)\pm2|> e^{- 10\Lambda q_k}, \ \ \forall k\ge K_2(E)
\end{equation}

\end{lemma}
With these two lemmas, we will first have the

\noindent
\textbf{Proof of Theorem \ref{thmtrace}:}

Follows immediately by combining Lemma \ref{lem1} with
$\tau=10\Lambda$ and  Lemma \ref{lem2}.
\qed

\begin{rmk}
It may be interesting to compare Theorem \ref{thmtrace} with the
technique Last used in his proof of zero Hausdorff dimensionality of
the spectral measures of supercritical Liouville almost Mathieu
operators \cite{l1}.  An important step there was using Schnol's theorem to show
that eventually spectrally almost every energy is in the union of the
spectral bands of the periodic approximants  {\bf enlarged} by a
factor of $q_k^2.$ Here we show that  spectrally almost every energy is in the union of the
{\bf shrinked} spectral bands of the periodic approximants, a much
more delicate statement, technically, thus with more powerful consequences.
\end{rmk}

\subsection{The hyperbolic case: Proof of Lemma \ref{lem1}}\label{sec2.3}
We are going to show that if $q$ is an `approximate' period  as in (\ref{beta-exp-peri}) with certain $\beta,\epsilon>0$ and satisfies
\begin{equation}
    |{\rm Trace} \,A_q\big(E\big)|\ge2+e^{- \tau q}
\end{equation}
then the trace of the transfer matrix at the scale $e^{\tau q/2 }$ will
be very large and any generalized eigenfunction of $Hu=Eu$ will be
bounded from below at the scale $e^{\tau q/2 }$. If this happens for
infinitely many $q$, then any generalized eigenfunction will have al
least a larger than $1/2$ power law growth (in index)
along some fixed subsequence.  By the extended Schnol's Theorem, such $E$ must belong to a set of spectral measure zero.

\begin{claim}\label{clmGrow}
Suppose $q\to\infty$ satisfy $|{\rm Trace} A_q\big(E\big)|\ge2+e^{- \tau q}$
and
\begin{equation}\label{qbeta}
     \max_{1\le j \le q, |k|\le e^{\epsilon\beta  q}/q}|V(j+kq)-V(j+(k\pm1)
q)|\le e^{-\beta q},\ \epsilon>0 .
\end{equation}
Assume further that
 \begin{equation}\label{betaepsilon}
    \beta>(3+1/\epsilon)\tau+(7+1/\epsilon)\Lambda
 \end{equation}
 Then there is $x_q^i\in\Z,i=1,...,4,$ independent of $E$, such that
 $|x_q^i|\to \infty$ as $q\to\infty$ and  for any $|u_0|^2+|u_1|^2=1$,
 $\max_{i=1,...,4}|u^E_{x_q^i}|>1/16e^q,$
 where $u^E_{n}$ is a solution with boundary values  $(u_0,u_1)$.
\end{claim}


\begin{lemma}[Extended Schnol's Theorem]\label{schnol}
Fix any $y>1/2$. For any sequence $|x_k|\to \infty$(where the sequence
is independent of $E$), for spectrally a.e. $E$, there is a generalized eigenvector $u^E$ of $Hu=Eu$, such that
$$|u^E_{x_k}|<C(1+|k|)^{y}$$
\end{lemma}
We can now finish the proof of Lemma \ref{lem1}.

\noindent \textbf{Proof of Lemma \ref{lem1}:}

Let $q_k$ be given as in (\ref{beta-exp-peri}). Combining Claim
\ref{clmGrow} with Lemma \ref{schnol} with
$\{x_k\}=\cup_{i=1}^{i=4}\{x_{q_k}^i\}$ we obtain that the set of $E$
such that there are infinitely many $q_{k_j}$ with $ |{\rm Trace}
A_{q_{k_j}}\big(E,\alpha\big)|\ge2+e^{- \tau q_{k_j}}$ has spectral
measure zero.
\qed

Claim \ref{clmGrow} is based on the following results. First, we need
to estimate the norm of the conjugation matrix for any hyperbolic $SL(2,\R)$ matrix w.r.t. the distance between its trace and $2$:
\begin{lemma}\label{lemG} Suppose $G\in SL(2,\R)$ with
  $2<|{\rm Trace}G|\le 6.$ The invertible matrix $B$ such that
\begin{equation}\label{G0}
    G=B\left(
                \begin{array}{cc}
                  \rho & 0 \\
                  0 & \rho^{-1} \\
                \end{array}
              \right)B^{-1}
\end{equation}
where $\rho^{\pm1}$ are the two conjugate real eigenvalues of $G$ with $|{\rm det}B|=1$ satisfies
\begin{equation}\label{B0}
    \|B\|=\|B^{-1}\|<\frac{\sqrt{\|G\|}}{\sqrt{|{\rm Trace}G|-2}}
\end{equation}
If $|{\rm Trace}G|> 6$,  then $\|B\|\le\frac{2\sqrt{\|G\|}}{\sqrt{|{\rm Trace}G|-2}}$.
\end{lemma}
The proof is based on a direct computation of the conjugate
matrices. For the sake of completeness, we present it in Appendix A.1.\\

Fix $\tau>0$ and apply Lemma \ref{lemG} to  $A_q$ satisfying
\begin{equation}\label{Etrace}
|{\rm Trace}\, A_q|>2+e^{-\tau q}.
\end{equation}We then have
\begin{claim}\label{clm1} For large $q,$
\begin{equation}\label{Aq}
    A_q=B\left(
                \begin{array}{cc}
                  \rho & 0 \\
                  0 & \rho^{-1} \\
                \end{array}
              \right)B^{-1}
\end{equation}
where $\rho^{\pm1}$ are the two conjugate real eigenvalues of $A_q$ with $\rho>1$ and $B$ satisfies $|{\rm det}B|=1$ and
\begin{equation}\label{B}
    \|B\|=\|B^{-1}\|<e^{(\tau/2+\Lambda+\Lambda/200)q}
\end{equation}
\end{claim}
Second, we need to use the almost periodicity (\ref{qbeta}) of the
potential to obtain approximation statements for the transfer matrices. Set
\begin{equation}\label{N}
    N=[e^{ (\tau/2+\Lambda/100) q}].
\end{equation}
Under the assumption (\ref{qbeta}) and (\ref{betaepsilon}) on $V$ as in Claim \ref{clmGrow}, we have for $q$ large enough (the largeness depend on $\Lambda$ and $n_0$),
\begin{claim}\label{clmNqNq} Under the conditions of Claim \ref{clmGrow},
\begin{equation}\label{NqNq}
    \|A_{Nq}-A_q^N\|\le 2e^{-\Lambda q}|\rho|^N\le 2e^{-\Lambda q} |{\rm Trace}A^N_q|
\end{equation}
and
\begin{equation}\label{-Nq-Nq}
    \|\big[A_{Nq}\big]^{-1}-A_{-Nq}\|\le 4e^{-\Lambda q} |\rho|^N \le 4e^{-\Lambda q} |{\rm Trace}A^N_q|
\end{equation}
\end{claim}
Let us now finish the proof of Claim \ref{clmGrow}.

\noindent \textbf{Proof of Claim \ref{clmGrow}}:
Decomposing $A_q$ as in (\ref{Aq}), we obtain
$|\rho|>
 1+e^{- \tau q/2}$. Obviously, $ |{\rm Trace}A_q^N| \ge  |\rho|^N$. By (\ref{N}), $N>2e^{ \tau q/2}\cdot q$, thus
$$|{\rm Trace}A_q^N|\ge (1+e^{- \tau q/2})^{2e^{ \tau q/2}\cdot q}
   \ge e^q.$$
Assume $q$ is large enough so that $2e^{-\Lambda q}\le 1/10$. By (\ref{NqNq}), we have,
\begin{equation}\label{trace1}
    |{\rm Trace}A_{Nq}|>(1-2e^{-\Lambda q})|{\rm Trace}A_q^N|\ge \frac{9}{10}e^q.
\end{equation}
Now consider  solution $u$ of $Hu=Eu$ with normalized initial value
$$X=\left(
  \begin{array}{c}
    u_1 \\
    u_0 \\
  \end{array}
\right), \ \ \|X\|=1$$
Then by (\ref{ref-sln}):
\begin{equation}\label{uq}
   A_{N q}\cdot X=\left(
  \begin{array}{c}
    u_{Nq+1} \\
    u_{Nq} \\
  \end{array}
\right), \ \ \ A_{-Nq}\cdot X=\left(
  \begin{array}{c}
    u_{-Nq+1} \\
    u_{-Nq} \\
  \end{array}
\right).
\end{equation}
By the Cayley-Hamilton theorem  combined with (\ref{-Nq-Nq}) and (\ref{trace1}), we have
\begin{eqnarray*}
\frac{9}{10}|{\rm Trace}A_q^N|\cdot \|X\|&\le&  \|{\rm Trace}A_{Nq}X\| \\
&=& \|A_{N q}\cdot X+\big[A_{Nq}\big]^{-1}\cdot X\| \\
 &\le& \|A_{N q}\cdot X\|+\|A_{-Nq}\cdot X\|+\frac{2}{10}|{\rm Trace}A_q^N|\cdot \|X\|
\end{eqnarray*}
Then
$$\|A_{N q}\cdot X\|+\|A_{-Nq}\cdot X\|\ge\frac{7}{10}|{\rm Trace}A_q^N|\cdot \|X\|\ge \frac{1}{2}|{\rm Trace}A_q^N|$$
which is equivalent to
$$\max\big\{\Big\|\left(
  \begin{array}{c}
    u_{Nq+1} \\
    u_{Nq} \\
  \end{array}
\right)\Big\|,\  \Big\|\left(
  \begin{array}{c}
    u_{-Nq+1} \\
    u_{-Nq} \\
  \end{array}
\right)\Big\|\big\}\ge 1/4 |{\rm Trace}A_q^N| .$$
Therefore
$$\max\Big\{\big|u_{N q+1}\big|,
   \big| u_{N q}\big|,
  \big|u_{-N q+1}\big|,
    \big|u_{-N q} \big|\Big\}\ge 1/16e^q.$$

Let $x_q^i=(-1)^iNq+1-[i/3],i=1,...,4.$ 
      Then for every $q$ and one of $i=1,...,4,\;$$|u_{x_q^i}|>1/16e^q$. \qed

It now remains  to prove (\ref{NqNq}),(\ref{-Nq-Nq}) in Claim \ref{clmNqNq}. 
Set
\begin{equation}\label{Delta}
    \Delta_i=A_{q}(T^{(i-1)q}\theta,E)-A_{q}(\theta,E), \ \ i=-N+1,\cdots,N.
\end{equation}
\begin{claim}\label{expTransfer}
Suppose (\ref{Lambda}) holds 
for $n\ge n_0$ and is uniform in $E\in S
$. Fix $E\in S,\theta\in\Omega$. If $V_\theta$ satisfies (\ref{qbeta}) with $\epsilon>0$, then there is a constant $C_{n_0}$ (depends only on $n_0$ and upper bound of $\|V\|_{\infty}$), such that
\begin{equation}\label{estDelta}
 \|\Delta_i(\theta,E)\| \le |i-1|qC_{n_0}e^{(\Lambda-\beta)q},\ \  |i|=1,\cdots,[e^{\epsilon \beta q}/q]
\end{equation}

\end{claim}
\proof
The proof is quite standard. Suppose $1\le i\le[e^{\epsilon \beta
  q}/q]$. Then for $|k|<i$, $|k|q<e^{\epsilon \beta q}$. Since
$V_{T^{kq}\theta}(n)=V_\theta(n+kq),$ (\ref{qbeta}) implies that for $|k| < i$ the following holds:
$$|V_{T^{kq}\theta}(j)-V_{T^{(k+1)q}\theta}(j)|\le e^{-\beta q},\ 1\le j\le q$$
which implies $$\|A(T^{kq+j}\theta)-A(T^{(k+1)q+j}\theta)\|\le e^{-\beta q},\ 1\le j\le q, \ |k| < i $$
By a standard telescoping argument, for any $\theta'=T^{kq}\theta, |k|<i$,
$$\|A_q(T^q\theta')-A_q(\theta')\|\le qC_V^{n_0}e^{(\Lambda-\beta) q}
=q C_{n_0}e^{(\Lambda-\beta)q}$$
where $C_V$ is such that $\|A(\theta',E)\|\le C_V,\forall \theta',E$.
In the above estimate, if 
$n>n_0$, we use the bound (\ref{Lambda}). When $n\le n_0$, we use the trivial bound $\|A_n\|\le C_V^{n_0}$,
We have
$$A_{iq}\big(\theta,E\big)=A_q\big(T^{(i-1)q}\theta\big)\cdots A_q\big(T^q\theta\big)A_q\big(\theta\big).$$
Therefore for $1\le i\le[e^{\epsilon \beta q}/q]$
$$
  \|\Delta_i\| \le \sum^{i-1}_{k=1}\|A_q\big(T^{kq}\theta\big)-A_q\big(T^{(k-1)q}\theta\big))\|
   \le (i-1)q C_{n_0}e^{(\Lambda-\beta)q}
$$
Since (\ref{qbeta}) is symmetric w.r.t. $T\to T^{-1}$, (\ref{estDelta}) for $i\le 0$ follows by takeing $T'=T^{-1}$. \qed

\noindent \textbf{Proof 
  of Claim \ref{clmNqNq}:}
Write  for any $i$, $A_q^i=B^{-1}\left(
                \begin{array}{cc}
                  \rho^i & 0 \\
                  0 & \rho^{-i} \\
                \end{array}
              \right)B$
so $\|A_q^i\|\le \|B\|^2|\rho|^i$. Set
$G(\theta)=\frac{1}{\rho}A_q(\theta)$ and
$G_j=G(T^{(j-1)q}\theta)$. By (\ref{B}), we have
$\|G^i\|\le \|B\|^2\le e^{(\tau+2\Lambda+\Lambda/100)q}$. Under the
assumption (\ref{betaepsilon}), we have $\tau/2+\Lambda/100<\epsilon
\beta$ so by (\ref{N}), $Nq<e^{\epsilon \beta q}$. Then Claim
\ref{expTransfer} implies, for $j=-N,\cdots, N$ and large $q,$ that
$$ \|G_j-G\|= \frac{1}{\rho}\|\Delta_j\| \le N qC_{n_0}e^{(\Lambda-\beta)q} \le e^{(-\beta+\tau/2+\Lambda+\Lambda/50)q} $$
Now we want to apply Lemma \ref{GN} to these $G_j$, with $M=e^{(\tau+2\Lambda+\Lambda/100)q}$ and $\delta=e^{(-\beta+\tau/2+\Lambda+\Lambda/50)q}$.  Direct computation gives
$$NM^2\delta<
e^{(-\beta+3\tau+5\Lambda+\Lambda/20)q}$$
By (\ref{betaepsilon}) we have $\beta-(3\tau+5\Lambda+\Lambda/20)>\Lambda$. Therefore, for $q$ large enough (the largeness only depends on $n_0$), we have $NM\delta<NM^2\delta<e^{-\Lambda q}$.
Then Lemma \ref{GN} implies that
$$\|\coprod_{j=1}^NG_j-G^N\|\le 2NM^2\delta\le 2e^{-\Lambda q}$$
and
$$\|\coprod_{j=1}^NG_{-N+j}-G^N\|\le 2NM^2\delta\le 2e^{-\Lambda q}$$

Therefore
\begin{equation}\label{N1}
    \|A_{Nq}-A_q^N\|=|\rho|^N\cdot\|\coprod_{j=1}^NG_j-G^N\|\le
2e^{-\Lambda q} |{\rm Trace}A^N_q|.
\end{equation}
establishing (\ref{NqNq}) ,
and
\begin{equation}\label{N2}
    \|A_{Nq}(T^{-Nq}\theta)-A_q^N(\theta)\|=|\rho|^N\cdot\|\coprod_{j=1}^NG_{-N+j}-G^N\|\le
2e^{-\Lambda q} |{\rm Trace}A^N_q|.
\end{equation}
Since $$\|A_{-Nq}(\theta)-\big[A^{-1}_{q}(\theta)\big]^N\|=\|\big[A_{Nq}\big]^{-1}(T^{-Nq}\theta)-\big[A^N_{q}(\theta)\big]^{-1}\|=\|\big[A_{Nq}\big](T^{-Nq}\theta)-\big[A^N_{q}(\theta)\big]\|$$
this implies
$$\|A_{-Nq}(\theta)-\big[A^{-1}_{q}(\theta)\big]^N\|\le 2e^{-\Lambda
  q} |{\rm Trace}A^N_q|.$$
Also,
\begin{equation}\label{-Nq1}
    \|\big[A_{Nq}\big]^{-1}-\big[A^{-1}_{q}\big]^N\|= \|\big[A_{Nq}\big]-\big[A_{q}\big]^N\|,
\end{equation}
therefore, by (\ref{N1}), we obtain  (\ref{-Nq-Nq}).
\qed

Lemma \ref{schnol}  is proved in the same way as the standard Schnol's
Lemma, however the statement in this form, while very useful, does not seem to be in the
literature (we learned it from S. Molchanov, see the Acknowledgement). For the sake of
completeness, we include a short proof in the Appendix.

\subsection{Energies with Trace close to $2$: Proof of Lemma \ref{lem2}}\label{sec2.4}
All throughout this section, we will assume again that all $q$ are large enough and satisfy (\ref{beta-exp-peri}) with certain $\beta,\epsilon>0$, i.e.,
\begin{equation}\label{qbeta2}
      \max_{1\le j \le q, |k|\le e^{\epsilon\beta  q}/q}|V(j+kq)-V(j+(k\pm1)
q)|\le e^{-\beta q},\ \epsilon>0.
\end{equation}
We are going to show that spectrally almost surely, there are only finitely many $q$ such that ${\rm Trace}A_{q}$ is close to $2$.

In fact, we are going to prove the following quantitative estimate on the measure of  energies where the trace of the associated transfer matrix is close to $2$.
\begin{lemma}\label{lem2'}
 Let $\Lambda$ be given by (\ref{Lambda}) on some set $S\subset \sigma(H)$. Let
\begin{equation}\label{Sqdef}
    S_q=\big\{E:\ 0<|{\rm Trace}A_q\pm2|<e^{- 10\Lambda q}\big\}.
\end{equation}
Assume (\ref{qbeta2}) holds and
\begin{equation}\label{tau3}
    \beta>(25+11/\epsilon)\Lambda.
\end{equation}
Then
\begin{equation}\label{Sqmes}
    \mu(S_q)<4q\cdot e^{-\Lambda q/15}<e^{-\Lambda q/20}
\end{equation}
where $\mu=\mu_S$ is the spectral measure restricted to $S$.
\end{lemma}

Once we have Lemma \ref{lem2'},  
Borel Cantelli lemma immediately  implies Lemma \ref{lem2}. So the
main task is to prove (\ref{Sqmes}).\\

In order to estimate the spectral measure of $S_q$, first we recall
the following results on the structure of $S_q.$ Let
$\mathcal{P}_n(\R)$ denote the polynomials over $\R$ of exact degree
$n$. Let  the class $\mathcal{P}_{n;n}(\R)$ be elements in $\mathcal{P}_n(\R)$ with $n$ {\em{distinct}} real zeros.
\begin{prop}[Theorem 6.1,\cite{jmx}] \label{jm}
Let $p \in \mathcal{P}_{n;n}(\R)$ with $y_1 < \dots < y_{n-1}$  the
local extrema of $p$. Let
\begin{equation} \label{eq_zeta}
\zeta(p) := \min_{1 \leq j \leq n-1} |p(y_{j})|
\end{equation} and $0 \leq  a < b$.  Then,
\begin{eqnarray}
|p^{-1}(a,b)|  \leq 2 \textrm{diam} (z(p-a)) \max\Big\{ \frac{b-a}{\zeta(p)+a}, \big(\frac{b-a}{\zeta(p)+a}\big)^\frac{1}{2} \Big\}
\end{eqnarray}

where $z(p)$ is the zero set of $p$ and $|\cdot |$ denotes the
Lebesgue measure.
\end{prop}
Fix any $\tau>0$. We apply Proposition \ref{jm} to polynomial ${\rm Trace}A_q(E)\in \mathcal{P}_{q;q}(\R),$ with $a=2$, $b=2+e^{- \tau q}$.
Clearly,  $\textrm{diam}(z({\rm Trace}A_q-2))$ is bounded from above by some
constant that only depends on $\|V\|_{\infty}$. We also have
$|\zeta({\rm Trace}A_q)|\ge 2$. Since $b-a<1$, we have $({\rm Trace}A_q)^{-1}(a,b)\le
C_V\sqrt{b-a}=C_V e^{- \tau q/2}$ where $C_V$ is some constant that
only depends on $\|V\|_{\infty}$. Since  $({\rm Trace}A_q)^{-1}(a,b)$ contains at most $q$ bands, setting $S_q=\big\{E:\ 2<{\rm Trace}A_q<2+e^{- \tau q}\big\}$, we have
\begin{equation}\label{SqIj}
    S_q=\bigcup_{j=1}^{q}{I}_j,\ \ \ |{I}_j|\le |S_q|\le C_V e^{- \tau q/2}.
\end{equation}The same analysis works for $(a,b)=(2-e^{- \tau q},2),(-2-e^{- \tau
  q},-2),(-2, -2+e^{- \tau q})$. Thus the structure (\ref{SqIj}) also
holds for the other three cases.

Denote by
\begin{equation}\label{Ij}
   \varepsilon^j _q= |{I}_j|<e^{(-\tau/2+\Lambda/200)q}.
\end{equation}

 If $ I_j\bigcap \Sigma\not=\emptyset,$ pick $E_j\in I_j\bigcap \Sigma$ where $\Sigma=\sigma(H)$ is the spectrum. Set $\widetilde{I}_j=(E_j-\varepsilon^j _q,E_j+\varepsilon^j _q)$. Then $I_j\subset \widetilde{I}_j$, so it is enough to estimate the spectral measure of  $\cup\widetilde{I}_j$.

Set $N_q=[e^{(\tau/2-\Lambda/200)q}]$.
For any $\varepsilon_q>0$, define $l_q=l(\varphi,\varepsilon_q,E),u^{\varphi},v^{\varphi}$ as in (\ref{l}). Write $l_q=[l_q]+l_q-[l_q]$, and $[l_q]=K_q\cdot q+r_q$, where $0\le r_q=[l_q]{\rm mod }\ q< q$
and $0\le l_q-[l_q]<1$. We need the following power law estimate, which is the key part to the proof of Lemma \ref{lem2'}.
\begin{claim}\label{clmuq}
Suppose $E\in S_q \bigcap \Sigma$ and $
 0<   \varepsilon_q<e^{(-\tau/2+\Lambda/200)q}$. Suppose
 (\ref{qbeta2}) holds. Assume that
 $\beta>(2+1/\epsilon)\tau+(5+1/\epsilon)\Lambda$ and $\tau\ge
 10\Lambda$. Then for every initial condition $\varphi$,
\begin{equation}\label{uvbd}
  \|u^{\varphi}\|_{l_q}^2\ge e^{\frac{1}{10}\Lambda q}
\end{equation}
\end{claim}
Combining (\ref{uvbd}) with the subordinacy theory, we are ready to
estimate the $m$-function and the spectral measure.\\

\noindent\textbf{Proof of Lemma \ref{lem2'}:} Take $\tau=10\Lambda$. Then $\beta>(25+11/\epsilon)\Lambda$ satisfies the requirement in Claim \ref{clmuq}.
Let $E_j\in I_j\bigcap \Sigma\subset S_q\bigcap \Sigma$ 
. For any $\varphi$, let $u^{\varphi,E_j}, v^{\varphi,E_j}$ be the
right half line solution associated with the energy $E_j$. According
to (\ref{Ij}), Claim \ref{clmuq} can be applied to all
$u^{\varphi,E_j}$.

We have for any $\varphi$,
$$\|u^{\varphi,E_j}\|_{l_q(j)}^2\ge e^{\frac{1}{10}\Lambda q}, \ \ j=1,\cdots,q$$
where $l_q(j)=l(\varphi,E_j,\varepsilon^j _q)$.

Then by the J-L inequality (\ref{jl}) and the definition of $l_q(j)$, we have
$$
\varepsilon^j _q\cdot |m_{\varphi}(E_j+ i\varepsilon^j _q)| <
  \frac{5+\sqrt{24}}{2\|u^{\varphi,E_j}\|_{l_q}\cdot\|v^{\varphi,E_j}\|_{l_q}}\cdot
  \frac{\|v^{\varphi,E_j}\|_{l_q}}{\|u^{\varphi,E_j}\|_{l_q}} < \frac{5+\sqrt{24}}{2}\cdot e^{-\Lambda q/10}
$$
Notice that the interval $I_j$ is independent of the boundary
condition $\varphi$, and so is $\varepsilon^j_q$. Therefore, we can
take the supremum w.r.t. $\varphi$ on both sides of the above
inequality. By  Proposition \ref{DKL}, we have
$$\varepsilon^j _q\cdot |M(E_j+ i\varepsilon^j _q)| 
\le \frac{5+\sqrt{24}}{2}\cdot e^{-\Lambda q/10} $$
On the other hand, by the definition of $M(z)$ in (\ref{M}), we have
$$\Im M(E+ i \varepsilon)\ge \frac{1}{2\varepsilon}\mu(E-\varepsilon,E+\varepsilon),\ \ E\in\R,\ \varepsilon>0$$
Therefore,
$$\mu(E_j-\varepsilon^j _q,E_j+\varepsilon^j _q)
\le 2\varepsilon^j _q\cdot |M(E_j+ i \varepsilon^j _q)|
\le (5+\sqrt{24}) e^{-\Lambda q/10}$$
which implies
$$\mu(I_j)\le\mu(\widetilde{I}_j)
\le  e^{-\Lambda q/15}$$
Since in (\ref{SqIj}) there are four cases for $S_q$ and each of them satisfies the previous estimates, the spectral measure of $S_q$ will be bounded by $4q e^{-\Lambda q/15}\le e^{-\Lambda q/20}$.  \qed

The proof of Claim \ref{clmuq} relies on the following estimates on the transfer matrices.
The first one is a 
formula for the power of a general $SL(2,\R)$ matrix. It is elementray
but turned out particularly useful and
will be an important part of our quantitative estimates in both
hyperbolic and nearly parabolic cases. As we did not find it in the
literature, we will provide a proof of it as well as of the next Lemma, in the Appendix.
\begin{lemma}\label{lemAkqFormula}
Suppose $A\in SL(2,\R)$ has eigenvalues $\rho^{\pm 1}$. For any $k\in\N$, if ${\rm Trace}A\neq2$, then
\begin{equation}\label{Akq}
    A^k=\frac{\rho^k-\rho^{-k}}{\rho-\rho^{-1}}\cdot\Big(
       A-\frac{{\rm Trace}A}{2}\cdot I\Big)+\frac{\rho^k+\rho^{-k}}{2}\cdot I
\end{equation}
Otherwise,
$A^k=k(A-I)+I$.
\end{lemma}
The key to the estimates in the nearly parabolic case is then the
following simple
\begin{lemma}\label{pho1}There are universal constants
  $1<C_1<\infty,c_1>1/3$ such that for $E\in S_q$ and $1\le k\le N_q$, we have
\begin{equation}\label{Akq1}
    c_1<\frac{\rho^k+\rho^{-k}}{2}<C_1\ ,\ \   c_1k<\frac{\rho^k-\rho^{-k}}{\rho-\rho^{-1}}<C_1k
\end{equation}
\end{lemma}
Second, since $A_q(\theta)$ is almost periodic  (with an exponential error),  the iteration of $A_q(\theta)$ along the orbit will be close to its power. The argument is similar to what we used in the proof of (\ref{NqNq}) in the previous part.
\begin{claim}\label{clmAq}
Fix $\theta\in\Omega$, $E\in S_q\bigcap \Sigma$ and $\tau>0$. Suppose
(\ref{qbeta2}) holds with $\beta>(2+1/\epsilon)\tau+(5+1/\epsilon)\Lambda$. Then for any $1\le k\le N_q$, we have
\begin{equation}\label{Akq-Akq}
    \|A_{kq}-A_q^k\|\le2e^{-\Lambda q}.
\end{equation}
\end{claim}
\proof
Set $\Delta_j=A_q(T^{j-1}\theta)-A_q(\theta)$. By the Claim
\ref{expTransfer}, $ \|\Delta_j\|\le jqCe^{(-\beta+\Lambda)q}, \
j<[e^{\epsilon \beta q}/q]$. Recall that
$N_q=[e^{(\tau/2-\Lambda/200)q}]$. The condition on $\beta$ guarantees
$N_q<[e^{\epsilon \beta q}/q]$, therefore we have
$\|\Delta_j\|\le e^{(-\beta+\tau/2+\Lambda+\Lambda/100)q}$ for all $j=1,\cdots, N_q.$  We need to check
the other requirements of Lemma \ref{GN}. According to Lemmas
\ref{lemAkqFormula}, \ref{pho1},
$$
  \|A_q^j\|<C_1j\|A_q-\frac{{\rm Trace}A_q}{2}\cdot I\|+C_1<3C_1N_q\|A_q\| <e^{(\tau/2+\Lambda+\Lambda/100)q}.$$

Now apply Lemma \ref{GN} to
the sequence
$A_q(\theta),\cdots,A_q(T^{j-1}\theta),\cdots,A_q(T^{k-1}\theta)$,
with $M=e^{(\tau/2+\Lambda+\Lambda/100)q}$ and
$\delta=e^{(-\beta+\tau/2+\Lambda+\Lambda/100)q}$.  We have
$N_qM^2\delta<
e^{(-\beta+2\tau+3\Lambda+\Lambda/40)q}.$
Since $\beta>(2+1/\epsilon)\tau+(5+1/\epsilon)\Lambda>2\tau+5\Lambda$
we have $\beta-(2\tau+3\Lambda+\Lambda/40)>\Lambda$. Therefore, for $q$ large enough,
$N_qM\delta<N_qM^2\delta<e^{-\Lambda q}$.Thus, by Lemma \ref{GN}, we have
$\|A_{kq}-A_q^k\|=\|\coprod_{j=1}^kA_q(T^{j-1}\theta)-A_q^k(\theta)\|\le 2e^{-\Lambda q}$.\qed

Now we are ready to finish the proof of the most technical part.\\

\noindent \textbf{Proof of Claim \ref{clmuq}:}
We first show the following lower bound for $K_q=[\frac{[\ell_q]}{q}]$:
\begin{equation}\label{Kq}
    K_q>e^{\Lambda q/6}>18C_1\cdot e^{\Lambda q/8}
\end{equation}
Actually, if $K_q\ge N_q$, 
(\ref{Kq}) is automatically satisfied since $\tau\ge 10\Lambda$.

Now assume $K_q<N_q$. For any $n\le [l_q]+1$, write $n=kq+r$, where $0\le k\le K_q, \ 0\le r\le q$.
Set $X_{\varphi}=\left(
                                                         \begin{array}{c}
                                                           \cos\varphi \\
                                                           -\sin\varphi \\
                                                         \end{array}
                                                       \right)
$. According to (\ref{Akq}), (\ref{Akq1}), we have
for any $\varphi$, $1\le k\le K_q<N_q,$
$$ \|A_q^k\cdot X_{\varphi}\|<C_1k\|A_q-\frac{{\rm Trace} A_q}{2}\cdot I\|+C_1< C_1k(\|A_q\|+3/2)+C_1$$
and by (\ref{Akq-Akq}),
$$\|A_{kq}\cdot X_{\varphi}\|\le\|A_q^k\cdot X_{\varphi}\|+
\|(A_{kq}-A_q^k)\cdot X_{\varphi}\|\le C_1k(\|A_q\|+3/2)+C_1+1$$
For $n_0< r\le q$, and for any $\theta'\in\Omega$,
$\|A_r(\theta')\|\le e^{\Lambda q}$. For $1\le r\le n_0$, we bound
$\|A_r(\theta')\|$ by $C^{n_0}$ as in the proof of Claim
\ref{expTransfer}. Therefore $\|A_r(\theta')\|\le e^{\Lambda q}$ for all $1\le r \le q$ with $q$ large. Thus,
\begin{eqnarray*}
  \|A_{kq+r}(\theta)\cdot X_{\varphi}\|
&\le& \|A_r(T^{kq}\theta)\|\cdot\|A_{kq}(\theta)\cdot X_{\varphi}\| \\
&\le & e^{\Lambda q}\Big(C_1k(\|A_q\|+3/2)+C_1+1\Big) \\
&\le& e^{\Lambda q}\cdot\Big(C_1k(e^{\Lambda q}+3/2)+C_1+1\Big)\\
&\le& k\cdot e^{(2\Lambda+\Lambda /200)q}
\end{eqnarray*}

Recall that $\left(
    \begin{array}{c}
      u^\varphi_{n+1} \\
      u^\varphi_{n} \\
    \end{array}
  \right)=A_n\cdot X_\varphi,
$ direct computation shows
\begin{eqnarray*}
   \|u^{\varphi}\|_{l_q}^2    &\le &  \sum_{n=1}^{[l_q]+1}\|A_n\cdot X_{\varphi}\|^2 \\
  &\le & \sum_{r=1}^{q}\|A_{r}\cdot X_{\varphi}\|^2
  +\sum_{k=1}^{K_q}\sum_{r=1}^{q}\|A_{kq+r}\cdot X_{\varphi}\|^2 \\
&\le& q\cdot e^{2\Lambda q} +
    \sum_{k=1}^{K_q}\sum_{r=1}^{q}k^2\cdot e^{(4\Lambda+\Lambda /100)q}\\
   &\le& q\cdot e^{2\Lambda q}+K_q^3\cdot q\cdot e^{(4\Lambda+\Lambda/100)q}\\
    &\le& K_q^3\cdot e^{(4\Lambda+\Lambda/20)q}
\end{eqnarray*}

Since $\varphi$ is arbitrary, and $\left(
        \begin{array}{c}
          v^{\varphi}_{n+1} \\
          v^{\varphi}_{n} \\
        \end{array}
      \right)=A_n\cdot X_{\varphi+\pi/2},$ $\|v^{\varphi}\|_{l_q}^2$ has the same upper bound. Therefore, $\|u^{\varphi}\|_{l_q}\|v^{\varphi}\|_{l_q}\le K_q^3\cdot e^{(4\Lambda+\Lambda/20)q}$. On the other hand, since $\varepsilon_q<e^{(-\tau/2+\Lambda/200)q}$, we have
\begin{equation}\label{lqlow}
\|u^{\varphi}\|_{l_q}\|v^{\varphi}\|_{l_q}=\frac{1}{2\varepsilon_q}\ge e^{(\tau/2-\Lambda/100)q}
\end{equation}
With $\tau\ge10\Lambda,$ this implies that $K_q^3
\ge e^{q\Lambda/2}$. Therefore,
\begin{equation}\label{k}
    K_q>e^{\Lambda q/6}>18C_1\cdot e^{\Lambda q/8}
\end{equation}
as claimed.

In order to get the lower bound on $\|u^{\varphi}\|_{l_q}^2$ in
(\ref{uvbd}), we need to consider  two cases.
\begin{description}
  \item[case I:] For $\varphi$ such that
\begin{equation}\label{case1}
    \|(A_{q}-\frac{{\rm Trace}A_q}{2}\cdot I)\cdot X_{\varphi}\|\ge e^{-\Lambda q/8},
\end{equation}
by (\ref{Akq}),(\ref{Akq1}), for any $1\le k\le18C_1\cdot e^{\Lambda
  q/8}\le N_q,$ we have
\begin{eqnarray*}
  \|A_{q}^k\cdot X_{\varphi}\|&=&
\|\frac{\rho^k-\rho^{-k}}{\rho-\rho^{-1}}\cdot\Big(
       A_q-\frac{{\rm Trace}A_q}{2}\cdot I\Big)X_{\varphi}+\frac{\rho^k+\rho^{-k}}{2}\cdot X_{\varphi}\|\\
   &\ge & \frac{\rho^k-\rho^{-k}}{\rho-\rho^{-1}}\cdot\|\Big(
       A_q-\frac{{\rm Trace}A_q}{2}\cdot I\Big)X_{\varphi}\|-\frac{\rho^k+\rho^{-k}}{2}\cdot \|X_{\varphi}\| \\
   &\ge& \frac{1}{3}k\cdot e^{-\Lambda q/8}-C_1
\end{eqnarray*}
Due to (\ref{Akq-Akq}), we have then
$$\|A_{kq}\cdot X_{\varphi}\|\ge
\|A_{q}^k\cdot X_{\varphi}\|-\|\Big(A_{kq}-A_{q}^k\Big)\cdot X_{\varphi}\|
\ge \frac{1}{3}k\cdot e^{-\Lambda q/8}-2C_1.$$
Therefore, for $9C_1\cdot e^{\Lambda q/8}\le k\le 18C_1\cdot e^{\Lambda q/8}$, we have
\begin{equation}\label{Akq4}
    \|A_{kq}\cdot X_{\varphi}\|\ge
C_1>1
\end{equation}
By (\ref{k}) and (\ref{Akq4}) we obtain
$$
   \|u^{\varphi}\|_{l_q}^2
   \ge \frac{1}{2} \sum_{n=1}^{[l_q]-1}\|A_n\cdot X_{\varphi}\|^2
   \ge \frac{1}{2}\ \ \sum_{k=[9C_1\cdot e^{\Lambda q/8}]+1}^{[18C_1\cdot e^{\Lambda q/8}]}\ \ \|A_{kq}\cdot X_{\varphi}\|^2
   \ge e^{\Lambda q/10}
$$

  \item[case II:] For $\varphi$ such that
\begin{equation}\label{case2}
    \|(A_{q}-\frac{{\rm Trace}A_q}{2} I)\cdot X_{\varphi}\|< e^{-\Lambda q/8},
\end{equation}
by (\ref{Akq}),(\ref{Akq1}), for any $1\le k\le N_q$ we get
\begin{eqnarray*}
  \|A_{q}^k\cdot X_{\varphi}\|&=&
\|\frac{\rho^k-\rho^{-k}}{\rho-\rho^{-1}}\cdot\Big(
       A_q-\frac{{\rm Trace}A_q}{2} I\Big)X_{\varphi}+\frac{\rho^k+\rho^{-k}}{2}\cdot X_{\varphi}\|\\
   &\ge & \frac{\rho^k+\rho^{-k}}{2}\cdot \|X_{\varphi}\|-\frac{\rho^k-\rho^{-k}}{\rho-\rho^{-1}}\cdot\|\Big(
       A_q-\frac{{\rm Trace}A_q}{2} I\Big)X_{\varphi}\| \\
   &\ge& 1/2-C_1k\cdot e^{-\Lambda  q/8}
\end{eqnarray*}
Combining with (\ref{Akq-Akq}), we have
$$\|A_{kq}\cdot X_{\varphi}\|\ge
\|A_{q}^k\cdot X_{\varphi}\|-\|(A_{kq}-A_{q}^k)\cdot X_{\varphi}\|
\ge \frac{1}{4}-C_1k\cdot e^{-\Lambda q/8}$$
Then for $1\le k\le \frac{1}{8C_1}\cdot e^{\Lambda q/8}\le K_q \le N_q$, we obtain
$\|A_{kq}\cdot X_{\varphi}\|
\ge 
\frac{1}{8}$.
This implies
$$
   \|u^{\varphi}\|_{l_q}^2
   \ge \frac{1}{2}\sum_{k=1}^{[\frac{1}{8C_1}\cdot e^{\Lambda q/8}]}\|A_{kq}\cdot X_{\varphi}\|^2 \\
   \ge e^{\Lambda  q/10}.
$$

\end{description}

\subsection{Proof of Lemmas \ref{lem4} and \ref{lem41}} \label{sectionlem4}
Assume that $|{\rm Trace}A_{q}|<2-e^{- \tau {q}}<2.$ Thus in the expression in Lemma \ref{lemAkqFormula}, $\rho=e^{ i\psi}, \psi\in(-\pi,\pi)$. We have for any $j$,
\begin{equation}\label{Akq5}
    A_q^j=\frac{\sin j\psi}{\sin \psi}\cdot\Big(
       A_q-\frac{{\rm Trace}A_q}{2}\cdot I\Big)+\frac{\cos j\psi}{2}\cdot I, \ \ \psi\in(-\pi,\pi)
\end{equation}
Then $|2\cos\psi|=|{\rm Trace}A_q|<2-e^{- \tau q}$ implies
$|\sin\psi|
>\sqrt{1-(1-\frac{1}{2}e^{- \tau q})^2}
>Ce^{- \tau q/2}$.
Therefore,
$$\|A_q^j\|\le C\cdot e^{ \tau q/2}\cdot\big(\|A_q\|+1\big)+1
$$
(here $q$ is large enough so that $\|A_q\|\le e^{\Lambda q}$.)
If $\tau=10\Lambda,$ we obtain
$$\|A_q^j\|\le e^{(6\Lambda+\Lambda/100)q}.
$$

Now let  $N=[e^{\xi q}]$. By the same argument as used for the
proof of (\ref{NqNq}) and (\ref{Akq-Akq}) (based on Lemma \ref{GN}), if $\beta>15\Lambda+(2+1/\epsilon)\xi$, then for any $j\le N^\xi$, $$\|A_q^j-A_{jq}\|<e^{(-\beta+13\Lambda+2\xi+\Lambda/20)q}<e^{-\Lambda q}.$$
As a consequence, we have $\|A_{jq}\|\le e^{(6\Lambda+\Lambda/50)q}$, and
$\|A_{jq+r}\|\le e^{(7\Lambda+\Lambda/50)q}$ for all $0\le r\le q$, $0\le j\le N^\xi$. Therefore
$$
  \sum_{n=1}^{N q}\|A_n(E)\|^2 \le\sum_{k=0}^{N^{\xi}}\sum_{r=1}^{q}\|A_{kq+r}(\theta,E)\|^2
   \le N q\cdot e^{(14\Lambda+\Lambda/25)q}
   \le e^{(\xi+15\Lambda)q}
$$\qed

{\bf Proof of Lemma \ref{lem41}}
Since $N>e^{(\xi-\Lambda/200) q}$ for $q$ large, then for any $\gamma<1$,
$$\frac{1}{(N q)^{2-\gamma}}\sum_{n=1}^{N q}\|A_n(E)\|^2<
e^{\big(-(1-\gamma)\xi+15\Lambda\big)q}
$$
If $\xi>\frac{16\Lambda}{1-\gamma},$ then $(1-\gamma)\xi-(15\Lambda)\ge 1/2\Lambda$. Therefore,
$$\frac{1}{(Nq)^{2-\gamma}}\sum_{n=1}^{Nq}\|A_n(E)\|^2\le e^{-1/2\Lambda q}\le1$$ \qed
\section{Spectral Singularity}\label{secSing}
\subsection{Power law estimates and proof of Theorem \ref{thmSingular}}
Throughout this section, our potential will be given by $V_\theta(n)=V(\theta+n\alpha),n\in\Z,$ where $V(\theta)$ is a real analytic function defined on the torus with analytic extension to the strip $\{z:|Imz|<\rho\}$.

According to Lemma \ref{lemMmm}, it is enough to find a $\varphi$ such
that both $m_\varphi$ and $\widetilde{m}_{\pi/2-\varphi}$ are $\gamma$-spectral singular. The main technical tool to estimate $m$-function is the subordinacy theory Lemma \ref{lemJL}.  We also need one more general statement about the existence of generalized eigenfunctions with sub-linear growth in its $l$-norm (see \cite{ls99}). That is, for $\mu_\theta$-a.e. $E$, there exists $\varphi\in (-\pi/2,\pi/2]$ such that $u^{\varphi}$ and $u^{\varphi,-}$ both obey
\begin{equation}\label{upperLS99}
    \limsup_{l\to\infty}\frac{\|u\|_l}{l^{1/2}\log l}<\infty
\end{equation}
This inequality provides us an upper bound for the $l$-norm of the
solution. To apply subordinacy theory, we also need a lower bound for the
$l$-norm. It will be derived from the following lower bounds on the norm of the transfer matrices. Denote
\begin{equation}\label{defAn-}
    \widetilde{A}_n(\theta,E,\alpha)=A_n(\theta-\alpha,E,-\alpha)
\end{equation}
We have
\begin{lemma}\label{lempoly}
Fix $\alpha\in\R\backslash\Q$ with $\beta=\beta(\alpha)<\infty$. Assume that $L(E)\ge a>0,E\in S$. There is $c=c(a,\rho)>0$ such that for $l>l(E,\beta,\rho)$, and any $\theta\in\T$, the following hold:
\begin{equation}\label{sumAk}
    \sum_{k=1}^{l} \|A_k(\theta,E,\alpha)\|^2\ge l^{1+\frac{2c}{\beta}}
\end{equation}
and
\begin{equation}\label{sumAk-}
    \sum_{k=1}^{l} \|\widetilde{A}_k(\theta,E,\alpha)\|^2\ge l^{1+\frac{2c}{\beta}}
\end{equation}
\end{lemma}

\noindent \textbf{Proof of Theorem \ref{thmSingular}:}
For any $\varphi$, we have
\begin{equation}\label{uvAk}
    \|u^{\varphi}\|_l^2+\|v^{\varphi}\|_l^2\ge \frac{1}{2}\sum_{k=1}^{l} \|A_k(\theta)\|^2
\end{equation}
and
\begin{equation}\label{uvAk-}
    \|u^{\varphi,-}\|_l^2+\|v^{\varphi,-}\|_l^2\ge \frac{1}{2}\sum_{k=1}^{l} \|\widetilde{A}_k(\theta)\|^2
\end{equation}
Therefore, a direct consequence of (\ref{sumAk}) is the power law estimate for the left hand side of (\ref{uvAk}), i.e.,
$\|u^{\varphi}\|_l^2+\|v^{\varphi}\|_l^2\ge l^{1+\frac{2c}{\beta}}$ for $l$ large.

On the other hand, according to (\ref{upperLS99}), for $\mu_\theta$-a.e. $E$, there exist $\varphi(E)$ and $C=C(E)<\infty$, such that for large $l$,
\begin{equation}\label{u}
    \|u^{\varphi}\|_l\le Cl^{1/2}\log l,\ \ \|u^{\varphi,-}\|_l\le Cl^{1/2}\log l
\end{equation}
Let us consider the right half line estimates for $u^\varphi,m_\varphi$ first. From (\ref{uvAk}) and (\ref{u}), we have
$$\|v^{\varphi}\|^2_l\ge l^{1+\frac{2c}{\beta}}-Cl(\log l)^2
$$
and then
\begin{equation}\label{v}
    \|v^{\varphi}\|_l\ge l^{1/2+c/{\beta}}
\end{equation}
provided $\beta<\infty$ and $l>l(\beta,E,\rho)$.\\

Applying subordinacy theory (\ref{jl}) to (\ref{u}),(\ref{v}), one has for any $\gamma\in(0,1)$, any $\varepsilon>0$
\begin{eqnarray*}
  \varepsilon^{1-\gamma}|m_\varphi(E+ i\varepsilon)| &\ge & \frac{1}{\Big(2 \|u^{\varphi}\|_{l(\varepsilon)}\|v^{\varphi}\|_{l(\varepsilon)}\Big)^{1-\gamma}}\cdot (5-\sqrt{24})\frac{\|v^{\varphi}\|_{l(\varepsilon)}}{\|u^{\varphi}\|_{l(\varepsilon)}}   \\
   &\ge& c_\gamma\frac{\|v^{\varphi}\|^{\gamma}_{l}}{\|u^{\varphi}\|^{2-\gamma}_{l}}  \\
   &\ge& c_\gamma l^{(1+c/{\beta})\gamma-1}\cdot \log^{-2}l
\end{eqnarray*}
where $c_\gamma>0$ may denote different constants that only depend on $\gamma$.
Set $\gamma_0=\gamma_0(\beta)=\frac{1}{1+c/\beta}<1,$ since $\beta<\infty$. We have for any $\gamma>\gamma_0$,
$$\varepsilon^{1-\gamma}|m_\varphi(E+ i\varepsilon)|\ge c_\gamma l^{\gamma/\gamma_0-1}\cdot \log^{-2}l\to\infty$$
as $\varepsilon\to 0$.

Using (\ref{uvAk-}) and (\ref{u}), the same argument works for $u^{\varphi,-},v^{\varphi,-}$ and $m^-_\varphi$.
Therefore, for spectrally a.e. $E$
$\liminf_{\varepsilon\to0}
\varepsilon^{1-\gamma}|m_\varphi(E+ i\varepsilon)|=\infty$ and
$\liminf_{\varepsilon\to0}
\varepsilon^{1-\gamma}|m^-_\varphi(E+ i\varepsilon)|=\infty$. According
to Lemma \ref{lemMmm},  $\mu$ is $\gamma$-spectral singular for any
$\gamma>\gamma_0$. The conclusion for the spectral dimension follows
 from the definition directly. \qed

The proof of Lemma \ref{lempoly} depends on the following lemmas about
the localization density of the half line solution. The key
observation is that in the regime of positive Lyapunov exponents we can guarantee transfer matrix growth at scale
$q_n$ somewhere within any interval of length $~q_n,$ giving a
contribution to (\ref{sumAk}).
\begin{lemma}\label{lemdegree}
Assume that $L(E)\ge a>0, \alpha\in\R\backslash\Q.$ There are $c_2=c_2(a,\rho)>0$ and a
positive integer $d=d(\rho)$ such that for $E\in S$  and
 $n>n(E,\rho)$, there exists an interval $\Delta_n$ such that
\begin{equation}\label{mea}
    Leb(\Delta_n)\ge \frac{c_2}{4dn}
\end{equation}
and for any $\theta\in \Delta_n$, we have
\footnote{
We denote by $\|\cdot\|_{HS}$ the Hilbert-Smith norm of a $SL(2,\R)$ matrix $$A=\left(
                                                                                                                   \begin{array}{cc}
                                                                                                                     a & b \\
                                                                                                                     c & d \\
                                                                                                                   \end{array}
                                                                                                                 \right)
,\ \ \|A\|_{HS}=\sqrt{a^2+b^2+c^2+d^2}$$
}
\begin{equation}\label{norm}
   \|A_n(\theta,E,\alpha)\|^2_{HS}> e^{n L(E)/8}
\end{equation}

\end{lemma}
In the following, we will use  $\|\cdot\|$  for the HS norm $\|\cdot\|_{HS}$. Now let $c_2$ and $d$ be given as in Lemma \ref{lemdegree}. Denote
\begin{equation}\label{kn}
    k_n=[\frac{c_2q_n}{4d}]-1
\end{equation}
where, as before $q_n$ are the denominators of the continued fraction
approximants to $\alpha.$
Based on Lemma \ref{lemdegree}, one can show that
\begin{lemma}\label{lemexp}
Fix $E\in S$  and  $\alpha\in\R\backslash\Q$. Let $k_n$ be given as in
(\ref{kn}). Suppose $q_n$ is large enough so
 that (\ref{mea}) holds for $\Delta_{k_n}$. Then for any $\theta$, and  any $N\in\N$, there is $j_N(\theta)\in[2Nq_n,2(N+1)q_n)$ such that
\begin{equation}\label{jn}
    \|A_{j_N}(\theta,E,\alpha)\|> e^{k_n L(E)/16}
\end{equation}

\end{lemma}
We first use Lemma \ref{lemdegree} and Lemma \ref{lemexp} to finish the proof of Lemma \ref{lempoly}. The proofs of these two lemmas are left to next section. \\

\noindent \textbf{Proof of Lemma \ref{lempoly}:} For $l$ sufficiently large, there is $q_n$ such that, $l\in[2q_n,2q_{n+1}).$ Write $l$ as
$$l=2Nq_n+r,$$
where $0\le r<2q_n$, $1\le N<\frac{q_{n+1}}{q_n}$. Suppose $q_n$ is large enough so that (\ref{mea}) holds for $\Delta_{k_n}$. Then Lemma \ref{lemexp} is applicable. Fix $\theta$. Consider $A_n(\theta,E,\alpha)$ first. Let $j_n(\theta)\in[2nq_n,2(n+1)q_n),n=0,1,\cdots,N,$ be given as in (\ref{jn}).
Direct computation shows that
\begin{eqnarray*}
  \sum_{k=1}^{l} \|A_k(\theta)\|^2 &\ge & \|A_{j_0}(\theta)\|^2+\|A_{j_1}(\theta)\|^2+\cdots+\|A_{j_{N-1}}(\theta)\|^2 \\
   &\ge&Ne^{k_n L(E)/16}
\end{eqnarray*}
We have $l=2Nq_n+r<4Nq_n$, i.e., $N>l/(4q_n)$. (\ref{kn}) implies
$c_2q_n/(5d)< k_n<c_2q_n/(4d)$ for $q_n$ large, so we have
$$\sum_{k=1}^{l} \|A_k(\theta)\|^2\ge \frac{l}{4q_n}e^{k_n L(E)/16}\ge 
\frac{l}{4q_n}e^{16c q_n} $$
where $c=c(c_2,d,a)$. Then for sufficiently large $l$, we have
$$\sum_{k=1}^{l} \|A_k(\theta)\|^2>l e^{8c q_n}$$
We also assume $l$ large enough (meaning $q_n$ large enough), so that
$\frac{\log q_{n+1}}{q_n} <2\beta$, i.e.,
$e^{q_n}>q_{n+1}^{\frac{1}{2\beta}}$.
Then
$$\sum_{k=1}^{l} \|A_k(\theta)\|^2\ge l\cdot q_{n+1}^{\frac{4c}{\beta}}\ge l\cdot (\frac{l}{2})^{\frac{4c}{\beta}}\ge l^{1+\frac{2c}{\beta}}. $$

For the same $\theta$, repeat the above procedure for
$A_n(\theta-\alpha,-\alpha,E)$. Notice
 that $\widetilde{A}_n(\theta,E,\alpha)=A_n(\theta-\alpha,E,-\alpha)$. Therefore, we have a sequence of positive integers $\widetilde{j}_N(\theta-\alpha)\in[2Nq_n,2(N+1)q_n)$ for any $N\in\N$ such that
\begin{equation}
    \|\widetilde{A}_{\widetilde{j}_N}(\theta,E,\alpha)\|> e^{k_n L(E)/16}
\end{equation}
The rest of the computations are exactly the same as for $A_n(\theta,E,\alpha)$. Notice that the constants $c_2$ and $d$ in Lemma \ref{lemdegree} are independent of the choice of $\alpha$ or $-\alpha$ and $\theta$. So $k_n$ and $c$ will be the same  for $A_n$ and $\widetilde{A}_n$
\qed


\subsection{Proof of the density lemmas}
\noindent \textbf{Proof of Lemma \ref{lemdegree}:}
Denote
\begin{equation}\label{fn}
    f_n(\theta)=\|A_n(\theta)\|^2_{HS}
\end{equation}
Obviously, $f_n(\theta)$ is a real analytic function with analytic
extension to the strip $\{z:|Imz|<\rho\}$. For bounded $S$ we have
\begin{equation}\label{fn2}
   \|f_n\|_{\rho}:=\sup_{|Imz|<\rho} \Big|f_n(z)\Big|<e^{C_1n}, \ \ E\in S
\end{equation}
where $C_1=C_1(S,\|V\|_\rho)$ can be taken uniform for all $E\in S$.
Expand $f_n(\theta)$ into its Fourier series on $\T$ as
\begin{equation}\label{fn3}
    f_n(\theta)=\sum_{k\in\Z}b_n(k)e^{2\pi ik\theta}
\end{equation}
where $b_n(k)$ is the $k$-th  Fourier coefficients of $f_n(\theta)$ so satisfies
\begin{equation}\label{bn}
    |b_n(k)|<\|f_n\|_{\rho}\cdot e^{-2\pi\rho|k|}, \ \forall k\in\Z
\end{equation}
We split $f_n(\theta)$ into two parts, for some positive integer $d$ which will be specified a little bit later
$$
  f_n(\theta)=g_n(\theta)+R_n(\theta),\
  g_n(\theta)=\sum_{|k|\le d\cdot n}b_n(k)e^{2\pi ik\theta},\ R_n(\theta)= \sum_{|k|>d\cdot n}b_n(k)e^{2\pi ik\theta}
$$
For any $\theta\in\T$
\begin{eqnarray*}
  |R_n(\theta)| \le  \sum_{|k|>d\cdot n}|b_n(k)| &\le &  \sum_{|k|>d\cdot n}\|f_n\|_{\rho}\cdot e^{-2\pi\rho|k|} \\
   &\le& \frac{2}{1-e^{-2\pi\rho}}e^{C_1n}e^{-2\pi\rho dn}
\end{eqnarray*}
Now pick
\begin{equation}\label{d}
    d=[\frac{C_1}{2\pi\rho}]+2
\end{equation}
With this choice of  $d$, we have $2\pi\rho d>C_1+1$, so for any $\theta\in\T$
\begin{equation}\label{rn}
    |R_n(\theta)|\le \frac{2}{1-e^{-2\pi\rho}}e^{-n}<1,\ \ n>n_0(\rho)
\end{equation}
Now we assume that the Lyapunov exponent $L(E)$ of $A(\theta,E)$ is positive. Denote
\begin{eqnarray*}
  \Theta_n^1 &=&\{\theta: f_n(\theta)> e^{n L(E)/8} \} \\
  \Theta_n^2 &=&\{\theta: g_n(\theta)> e^{n L(E)/4} \} \\
  \Theta_n^3 &=&\{\theta: f_n(\theta)> e^{n L(E)/2} \}
\end{eqnarray*}
According to (\ref{rn}), we see that if $f_n(\theta)> e^{n L(E)/2}$, then
$$g_n(\theta)>f_n(\theta)-|R_n(\theta)|>e^{n L(E)/2}-1>e^{n L(E)/4},\ \ n>n(E)$$
and if $g_n(\theta)> e^{n L(E)/4}$, then
$$f_n(\theta)>g_n(\theta)-|R_n(\theta)|>
e^{n L(E)/8},\ \ n>n(E)$$
Therefore, we have for large $n,$
\begin{equation}\label{theta}
   \Theta_n^3\subseteq\Theta_n^2\subseteq\Theta_n^1
\end{equation}
On the other hand,
\begin{eqnarray*}
  2nL(E) &\le& 2\int_{\T}\log\|A_n(\theta)\|_{HS}{\rm d}\theta=\int_{\T}\log f_n(\theta){\rm d}\theta \\
   &\le&  Leb(\Theta_n^3) \log \|f_n\|_{\rho}+ \big(1-Leb(\Theta_n^3)\big) \log e^{n L(E)/2} \\
   &\le&  Leb(\Theta_n^3) \cdot C_1n+ \big(1-Leb(\Theta_n^3)\big) \cdot n L(E)/2
\end{eqnarray*}
which implies $Leb(\Theta_n^3)\ge \frac{3L(E)}{2C_1-L(E)}$. Since
$L(E)\ge a>0,E\in S$, we have
\begin{equation}\label{ce}
    Leb(\Theta_n^3)\ge \frac{3a}{2C_1-a}:=c_2(a,\rho)>0
\end{equation}
Thus
\begin{equation}\label{ce2}
    Leb(\Theta_n^2)
\ge c_2(a,\rho)>0
\end{equation}
Since $g_n(\theta)$ is a trigonometric polynomial of degree
$2dn$, the set $\Theta_n^2$ consists of no more than $4dn$ intervals. Therefore, there exists a
segment, $\Delta_n\subset \Theta_n^2$, with $Leb(\Delta_n)>\frac{c_2}{4dn}$. Obviously, $\Delta_n$ is also contained in $\Theta_n^1$, i.e.,
for any $\theta\in \Delta_n$, $$\|A_n(\theta)\|^2_{HS}> e^{n L(E)/8}$$ and
$$Leb(\Delta_n)>\frac{c_2}{4dn},\ \ n>n(E,\rho) $$
where $d$ only depends on $\rho$ and is independent of $n$. \qed
The following standard lemma is proved e.g.  in \cite{jl2}
\begin{lemma}[Lemma 9, \cite{jl2}]\label{jl9}
Let $\Delta\subset [0,1]$ be an arbitrary segment.  If $|\Delta|>\frac{1}{q_n}$. Then, for any $\theta$; there exists a $j$ in
$\{0,1,\cdots,q_n+q_{n-1}-1\}$ such that $\theta+j\alpha\in\Delta$.
\end{lemma}

\noindent \textbf{Proof of Lemma \ref{lemexp}:}
The case $N=0$ is already covered by Lemma \ref{jl9}. The proof for the
case $N>0$ follows the same strategy.
Notice that (\ref{kn}) implies
$|\Delta_{k_n}|>\frac{c_2}{4dk_n}>\frac{1}{q_n}$ for large
$q_n$. Applying Lemma \ref{jl9} to $\theta+2Nq_n$, we have that there exists a $j$ in
$\{0,1,\cdots,q_n+q_{n-1}-1\}$ such that $\theta+2Nq_n\alpha+j\alpha\in\Delta_{k_n}$, i.e.,
$$\|A_{k_n}(\theta+2Nq_n\alpha+j\alpha)\|> e^{k_n L(E)/32}$$
Since
$$A_{2Nq_n+j+k_n}(\theta)=A_{k_n}(\theta+2Nq_n\alpha+j\alpha)A_{2Nq_n+j}(\theta)$$
and $A_i$ is unimodular, we have that either
$$\|A_{2Nq_n+j}(\theta)\|\ge e^{k_n L(E)/32}\ \ \textrm{or} \ \ \|A_{2Nq_n+j+k_n}(\theta)\|\ge e^{k_n L(E)/32}$$
Let $j_N$ be $2Nq_n+j$ or $ 2Nq_n+j+k_n$, so that $j_N$ satisfies (\ref{jn}).
Clearly,
$$2Nq_n\le 2Nq_n+j< 2Nq_n+j+k_n 
< 2Nq_n+2q_n$$
Therefore, $j_N\in[2Nq_n,2(N+1)q_n)$. \qed
\section{Sturmian Hamiltonian}\label{secSturm}

 Liu, Qu and Wen \cite{lw04, lqw14} studied the Hausdorff and upper box counting dimension of $\Sigma_{\lambda,\alpha}$ with general irrational frequencies. For any irrational $\alpha\in(0,1)$  with continued fraction expansion $[0;a_1,a_2,\cdots]$, define
\begin{equation}\label{K-u-l}
K_\ast(\alpha)=\liminf_{k\rightarrow\infty}
(\prod_{i=1}^k a_i)^{1/k}\ \text{ and }\  K^\ast(\alpha)=
\limsup_{k\rightarrow\infty}(\prod_{i=1}^k a_i)^{1/k}.
\end{equation}
Then (Theorem 1 \cite{lw04}, Theorem 1.1 \cite{lqw14})  for large coupling constant $\lambda$, $\dim_H \Sigma_{\alpha,\lambda} =1 $ iff $K_\ast(\alpha)= \infty$ and $\overline{\dim}_B \Sigma_{\alpha,\lambda} =1$ iff $K^\ast(\alpha)= \infty$.

The usual way to study Sturmian Hamiltonian is to decompose Sturmian potentials into canonical words, which obey recursive relations. Here we present an alternative approach to study spectral dimension properties of Sturmian Hamiltonian based on the technics we developed 
in Theorem \ref{thmConti}.\\

We will first prove Theorem \ref{thmSturm}. Set
\begin{equation}\label{Vsturm}
    V_\theta(n)=\lambda\chi_{[1-\alpha,1)}(\theta+n\alpha \ \textrm{mod}\ 1)
\end{equation}
It is well known that for Sturmian $H_\theta$, the restriction of Lyapunov exponent on the spectrum is
zero (see Theorem 1,
\cite{dl}). 
By the discussion  after (\ref{uniupper}) (see \cite{jmv2}) or else,
specifically for Sturmian potentials, by
\cite{lenz}, for arbitrarily small $\Lambda>0$ and $n\ge n_0(\Lambda)$, $\|A_n(\theta,E)\|\le e^{\Lambda n} $ uniformly in $\theta$ and $E\in\sigma(H_\theta)$. 
Here we will  apply Corollary \ref{zeroLyp} directly.

Let $q_k$ be the subsequence of denominators of the continued fraction
approximants of $\alpha$ such that $\|q_k\alpha\|<e^{-\beta q_k/2}$
. In order to apply Corollary \ref{zeroLyp}, it is enough to verify
that $V_{\theta}(n)$ given by (\ref{Vsturm}) is $\beta(\alpha)$-almost periodic for $\alpha$-Diophantine $\theta\in\T$. Fix  $\tau>1$. 
If $\theta$ is $\alpha$-Diophantine
there is $\gamma>0$ such that $\|\theta+m\alpha\|_{\R/\Z}\geq \frac{\gamma}{(|m|+1)^\tau}$ for any $m\in\Z$.
Then for $|m|\le q$
$$\textrm{dist}(\theta+m\alpha,\{\Z,1-\alpha+\Z\})\ge \min_{|m|\le q+1}\|\theta+m\alpha\|_{\R/\Z}.$$
Therefore,
$$\min_{|m|\le q}\textrm{dist}(\theta+m\alpha,\{\Z,1-\alpha+\Z\})
\ge \min_{|m|\le q+1}\frac{\gamma}{(|m|+1)^\tau} \ge \frac{\gamma}{(q+2)^\tau}.$$

 Let $N=[e^{\beta q/4}]$. Then for $|j|\le N$, $q>q_0(\gamma,\beta)$ and any $|m|\le q$, we have
$$\|jq\alpha\|\le |j|\cdot \|q\alpha\|\le e^{-\beta q/4}\le \frac{\gamma}{10(q+2)^\tau}\le
\frac{1}{10}\textrm{dist}(\theta+m\alpha,\{\Z,1-\alpha+\Z\})$$
Therefore,  for any $|m|\le q$ and
$|j|\le N$, $\theta+m\alpha \ \textrm{mod}\ 1$ and $\theta+m\alpha+jq\alpha \ \textrm{mod}\ 1$ belong to one of the same open intervals $\Big\{(0,1-\alpha),\ (1-\alpha,1)\Big\}$, which implies that
$$ \chi_{[1-\alpha,1)}(\theta+m\alpha \ \textrm{mod}\ 1)=\chi_{[1-\alpha,1)}(\theta+m\alpha+jq\alpha \ \textrm{mod}\ 1),
\ \ |m|\le q,|j|\le N$$
Therefore, for $0\le m \le q,$
$$V_\theta(m)=V_\theta(m+q)=\cdots=V_\theta(m+Nq)$$
which immediately implies $\beta(\alpha)$-almost periodicity for the sequence $q_k$ with $\epsilon=1/4$.

Since the set of $\alpha$-Diophantine $\theta$ has full Lebesgue
measure, the conclusion for the  density of states  follows directly from $\textrm{d}N=\E(\textrm{d}\mu_\theta)$.\\

Next we will construct $\alpha$ to prove Corollary \ref{corSturm}.
We will define inductively the continued fraction coefficients
$ a_n, n\geq 1,$ so $\alpha=[a_1,\cdots,a_{n},\cdots]$. Fix $\beta>0$. Start with some $n_0$ large. For $ 1\le i\le n_0$, set
$a_i=1.$  Set
$[a_1,\cdots,a_{n}]= p_{n}/q_{n}.$ Now, for $k=1,\cdots$ define $n_k=q_{n_0}+q_{n_1}+\cdots+q_{n_{k-1}}$
and
$a_{n}=\{\begin{array}{l}
         e^{\beta q_{n_k}};  \ n=n_k+1, \\
         1; \ \ \qquad  n_k+2\leq n\leq n_{k+1},
       \end{array}
$ for $\; k=0,1,\cdots$

Set $\alpha=[a_1,\cdots,a_{n},\cdots]$. It is  easy to check that
\begin{itemize}
  \item $$\beta+\frac{\log q_{n_k}}{q_{n_k}}=\frac{\log a_{n_k+1}q_{n_k}}{q_{n_k}}<\frac{\log q_{n_k+1}}{q_{n_k}}<
    \beta+\frac{\log 2q_{n_k}}{q_{n_k}}$$
  $$\Longrightarrow\frac{\log q_{n_k+1}}{q_{n_k}}\rightarrow \beta$$
  \item
  \begin{eqnarray*}
   (a_1a_2\cdots a_{n_k})^{1/{n_k}}&=&
  (a_{n_0+1}a_{n_1+1}\cdots a_{n_{k-1}+1})^{1/n_k} \\
     &= & (e^{\beta q_{n_0}}e^{\beta q_{n_1}}\cdots e^{\beta q_{n_{k-1}}})^{1/(q_{n_0}+q_{n_1}+\cdots+q_{n_{k-1}})} \\
     &=& e^\beta<\infty
  \end{eqnarray*}
\end{itemize}

Therefore, $\alpha$ constructed in the above way satisfies
$\beta(\alpha)>0$ while $K_\ast(\alpha)<\infty$. Then Corollary
\ref{corSturm} follows from \cite{lw04} and Theorem \ref{thmSturm}. \qed

On the other hand, if we take $\alpha=[0;1,2,3,\cdots,k,\cdots]$, then $K_\ast(\alpha)=\infty$ while $\beta(\alpha)=0$. By \cite{lw04}, for Sturmian Hamiltonian with frequencies $\alpha$ such that $K_\ast(\alpha)=\infty$, $\dim_H \Sigma_{\alpha,\lambda}=\dim_P\Sigma_{\alpha,\lambda} =1 $. 

\appendix
\section{Appendix}
\subsection{Proof of Lemma \ref{lemG}:}
Suppose that $u,v$ are the two normalized eigenvectors of $G$ such that
$$Gu=\rho u,\ Gv=\rho^{-1}v,\ \ \|u\|=\|v\|=1$$
Denote the angle between $u$ and $v$ by $\theta$. Without loss of generality we assume further that $|\theta|<\pi/2$. Set
$\widetilde{B}=(u,v),\ \ B=\frac{\widetilde{B}}{\sqrt{|det \widetilde{B}|}}$.
Obviously, $\|\widetilde{B}\|\le 1$, $|det B|=1$, and $det \widetilde{B}=\|u\|\cdot \|v\|\cdot \sin \theta$. Therefore,
$$\|B\|\le \frac{1}{\sqrt{|\sin \theta|}}$$
On the other hand,
$G(u-v)=\rho u-\rho^{-1}v$, which implies that
$$\rho -\rho^{-1}=\rho \|u\|-\rho^{-1}\|v\|\le \|\rho u-\rho^{-1}v\|=\|G(u-v)\|\le \|G\|\cdot \|u-v\|$$
By the law of cosines, $\|u-v\|
=2\sin\frac{\theta}{2}$. Then
$$2\sin\frac{\theta}{2}\ge \frac{\rho-\rho^{-1}}{\|G\|}= \frac{\sqrt{(|{\rm Trace}G|+2)(|{\rm Trace}G|-2)}}{\|G\|}$$
$|{\rm Trace}G|\le 6$ implies that $|{\rm Trace}G|+2\ge2(|{\rm Trace}G|-2)$, then $2\sin\frac{\theta}{2}\ge \frac{\sqrt2(|{\rm Trace}G|-2)}{\|G\|}$.

Therefore,
$$\sin \theta\ge 2\sin\frac{\theta}{2}\cdot \frac{1}{\sqrt2}
\ge \frac{|{\rm Trace}G|-2}{\|G\|}$$
and
$$\|B\|\le\frac{\sqrt{\|G\|}}{\sqrt{|{\rm Trace}G|-2}}$$
It is also easy to see that if $|{\rm Trace}G|> 6$, $\|B\|\le\frac{2\sqrt{\|G\|}}{\sqrt{|{\rm Trace}G|-2}}$.
 \qed

\subsection{Proof of Lemma \ref{lemAkqFormula} and Lemma \ref{pho1}}
\noindent \textbf{Proof of Lemma \ref{lemAkqFormula}:}
Suppose $A=\left(
                \begin{array}{cc}
                  a & b \\
                  c & d \\
                \end{array}
              \right)\in SL(2,\R)$ has eigenvalues $\rho$ and
              $\rho^{-1}.$

\noindent \textbf{Case I:} ${\rm Trace}A\neq2$. Obviously, $\rho\neq 1$ and
\begin{equation}\label{ApdAq}
    A=\left(
                \begin{array}{cc}
                  a & b \\
                  c & d \\
                \end{array}
              \right)=B\left(
                \begin{array}{cc}
                  \rho & 0 \\
                  0 & \rho^{-1} \\
                \end{array}
              \right)B^{-1}
\end{equation}
where $B$ is the conjugation matrix. Suppose $\rho\neq d$. We can pick the conjugation matrix as
\begin{equation}\label{ApdB}
    B=\left(
                \begin{array}{cc}
                  1 & \frac{b}{\rho^{-1}-a} \\
                   \frac{c}{\rho-d}  & 1 \\
                \end{array}
              \right),\ \ \ B^{-1}=\frac{\rho-d}{\rho-\rho^{-1}}\left(
                \begin{array}{cc}
                  1 & -\frac{b}{\rho^{-1}-a} \\
                   -\frac{c}{\rho-d}  & 1 \\
                \end{array}
              \right).
\end{equation}
If $\rho=d$, it is easy to see that $bc=0$. Without loss of generality, we assume $c=0,b\neq0$. We can pick the conjugation matrix as
\begin{equation}\label{ApdB2}
    B=\left(
                \begin{array}{cc}
                  1 & 1 \\
                   \frac{d-d^{-1}}{b}  & 0 \\
                \end{array}
              \right),\ \ \ B^{-1}=\frac{b}{d^{-1}-d}\left(
                \begin{array}{cc}
                  0 & -1 \\
                   -\frac{d-d^{-1}}{b}  & 1 \\
                \end{array}
              \right).
\end{equation}

Direct computation using (\ref{ApdAq}),(\ref{ApdB}),(\ref{ApdB2}) shows that for any $k\in\N$,
\begin{equation}\label{ApdAqk}
  A^k=\frac{\rho^k-\rho^{-k}}{\rho-\rho^{-1}}\cdot\Big(
       A-\frac{a+d}{2}\cdot I\Big)
+\frac{\rho^k+\rho^{-k}}{2}\cdot I
\end{equation}

\noindent \textbf{Case II:} ${\rm Trace}A=2$.
Also follows by a (simpler) direct computation, considering separately $a=1$ and $a\neq1$.
\qed

\noindent \textbf{Proof of Lemma \ref{pho1}:}
Now assume $E\in S_q$ and $1\le k\le N_q$. Apply (\ref{ApdAqk}) to $A_q(E)$. First, suppose $2<{\rm Trace}A_q(E)<2+e^{- \tau q}$.
Then
$$1<\rho=\frac{{\rm Trace}A_q(E)+\sqrt{({\rm Trace}A_q(E))^2-4}}{2}<\frac{2+e^{- \tau q}+\sqrt{(2+e^{- \tau q})^2-4}}{2}<1+e^{(-\tau/2+\Lambda/200)q}$$

There is a universal
constant $C$, such that for any $ 1\le k\le N_q<e^{(\tau/2-\Lambda/200)q}$,
 $$1<\rho^k<(1+e^{(-\tau/2+\Lambda/200)q})^{N_q}<C.$$
Therefore,
\begin{equation}
    1<\frac{\rho^k+\rho^{-k}}{2}<C.
\end{equation}
On the other hand,
$$\frac{\rho^k-\rho^{-k}}{\rho-\rho^{-1}}
   =\sum_{i=1}^{k}\rho^{k-2i+1}$$
therefore,
\begin{equation}
    k\leq\frac{\rho^k-\rho^{-k}}{\rho-\rho^{-1}}<C_1k
\end{equation}

Now assume $2-e^{- \tau q}<{\rm Trace}A_q(E)<2$. Then
$\rho=e^{i\psi}$ and
(\ref{ApdAqk}) can be expressed as
\begin{equation}
   A_q^k
  = \frac{\sin k\psi}{\sin \psi}\cdot\Big(
       A_q-\frac{{\rm Trace}A_q}{2}\cdot I\Big)+\frac{\cos k\psi}{2}\cdot I
\end{equation}
We have $1-\frac{1}{2}e^{- \tau q}<\cos \psi<1$.
Then $|\sin\psi|
<e^{- \tau q/2}$, and
$|\psi|<\frac{\pi}{2}|\sin\psi|<2e^{- \tau q/2}$. As in the hyperbolic case, we set $N_q=[e^{(\tau/2-\Lambda/200)q}]$.
For $k\le N_q$, $$|k\psi|<
2e^{-\Lambda q/200}$$
Then for $q$ large enough, we have $
  \frac{2}{\pi} |k\psi|\le  |\sin k\psi|\le |k\psi|\le \sqrt3/2
$. Therefore
$
  \frac{2}{\pi}k\le
\Big|\frac{\sin k\psi}{\sin \psi}\Big|<\frac{\pi}{2}k
$
and
$
    1\ge \cos k\psi>1/2
$.

Exactly the same argument works for the case $\big\{E:\ -2<{\rm Trace}A_q<-2+e^{- \tau q}\big\}$ and $\big\{E:\ -2-e^{- \tau q}<{\rm Trace}A_q<-2\big\}$. \qed


\subsection{Some estimates on matrix products}
\begin{lemma}\label{GN}
Suppose $G$ is a two by two matrix satisfying
\begin{equation}\label{G^j}
    \|G^j\|\le M<\infty, \ \ \ \textrm{for all}\ 0<j\le N,
\end{equation}
where $M\ge 1$ only depends on $N$. Let $G_j=G+\Delta_j$,
$j=1,\cdots,N,$ be a sequence of two by two matrices with
\begin{equation}\label{delta}
    \delta=\max_{1\le j\le N}\ \|\Delta_j\|.
\end{equation}
If
\begin{equation}\label{deltasmall}
    NM\delta<1/2,
\end{equation}
then for any $n\le N$
\begin{equation}\label{G-product}
    \|\coprod_{j=1}^nG_j-G^n\|\le 2NM^2\delta
\end{equation}
\end{lemma}

\noindent \textbf{Proof of Lemma \ref{GN}:}
Denote by $$D=\max_{1\le k_1,k_2\le N}\|\coprod_{j=k_1}^{k_2}G_j\|.$$
Then a simple perturbation argument, as in e.g.  \cite{sb}, one can show that $D\le M(\delta D N+1)$. Thus
$D\le \frac{M}{1-M\delta N}$.
Direct computation shows that for any $1\le n\le N$,
$$\coprod_{j=1}^nG_j-G^n=\sum_{k=0}^{n-1}\big(\coprod_{j=k+2}^{n}G_j\big)\Delta_{k+1} G^k$$
Therefore, $$\|\coprod_{j=1}^nG_j-G^n\|\le ND\delta M\le \frac{M^2\delta N}{1-M\delta N}$$
Clearly, if $M\delta N<1/2$, then $\|\coprod_{j=1}^nG_j-G^n\|\le2NM^2\delta.$ \qed


\subsection{Extended Schnol's Theorem (Lemma \ref{schnol})}
Let $y>1/2$ and  $x_k$be any sequence such that  $|x_k|\to\infty$ as $k\to\infty$. For a Borel set $B\in\R$, denote
$$\mu_{n,m}(B)=<\delta_n,\chi_B(H)\delta_m>$$
and
$$\rho(B)=\sum_{n}a_n(\mu_{n,n}(B)+\mu_{n+1,n+1}(B))$$
where
$$a_n=\left\{
    \begin{array}{ll}
     c(1+|k|)^{-2y}, & n=x_k \\
      c(1+|n|)^{-2y}, & else
    \end{array}
  \right.
$$ with $c > 0$ chosen so that $\sum_na_n=1/2$.

Then, $\rho$ is a Borel probability measure with $\rho(B)=0$ if and
only if $\mu(B) = 0$, i.e., $\rho$ and $\mu$ are mutually absolutely
continuous. By the Cauchy-Schwarz inequality,
$$|\mu_{n,m}(B)|\le \mu_{n,n}(B)^{\frac{1}{2}}\mu_{m,m}(B)^{\frac{1}{2}}.$$
therefore
$\mu_{n,m}$ is absolutely continuous
with respect to $\rho$. By the Radon-Nikodym Theorem, there exists a measurable
density
$$F_{n,m}(E)=\big[\frac{{\rm d}\mu_{n,m}}{{\rm d}\rho}\big](E), \ \ \rho-a.e. \ E $$ with
$$\mu_{n,m}(B)=\int\chi_B(E)F_{n,m}(E){\rm d}\rho(E).$$
Then for every bounded measurable
function $f$, we have that
$$<\delta_n,f(H)\delta_m>=\int f(E)F_{n,m}(E){\rm d}\rho(E)$$
In particular, if $g$ is compactly supported and bounded, we may set $f(E) = Eg(E)$
and have
\begin{eqnarray*}
  &&\int g(E)\Big(EF_{n,m}(E)\Big){\rm d}\rho(E)\\
   &=& <\delta_n,Hg(H)\delta_m> \\
   &=& <\delta_{n+1}+\delta_{n-1}+V_n\delta_{n},g(H)\delta_m> \\
   &=& \int g(E)F_{n+1,m}(E){\rm d}\rho(E)+\int g(E)F_{n-1,m}(E){\rm d}\rho(E)+\int g(E)V_nF_{n,m}(E){\rm d}\rho(E) \\
   &=& \int g(E)\Big[F_{n+1,m}(E)+F_{n-1,m}(E)+V_nF_{n,m}(E)\Big]{\rm d}\rho(E)
\end{eqnarray*}
For any fixed $m\in\Z$, let $u^E(n)=F_{n,m}(E).$ Thus we have for any $g$
$$\int g(E)\Big((H-E)u^E\Big)(n){\rm d}\rho(E)=0$$
i.e., $\{u^E(n)\}_{n\in\Z}$ is a generalized eigenfunction of $Hu=Eu$ for $\rho$ a.e. $E$. \\

On the other hand, 
let
$$B_n=\{E:\ F_{n,n}\ge\frac{1}{a_n}\}$$
Then $$\rho(B_n)=\sum_{k}a_k\mu_{k,k}(B_n)\ge a_n\mu_{n,n}(B_n)=a_n\int_{B_n}F_{n,n}(E){\rm d}\rho(E)$$
While
$$\int_{B_n}F_{n,n}(E){\rm d}\rho(E)\ge \frac{1}{a_n}\rho(B_n)$$
Therefore,
$$\int_{B_n}\Big(a_nF_{n,n}(E)-1\Big){\rm d}\rho(E)\leq 0$$
Therefore, $\rho(B_n)=0$, i.e., for $\rho$-a.e. $E$, $F_{n,n}(E)\le\frac{1}{a_n}$, thus
$$|F_{n,m}|\le a_n^{-\frac{1}{2}}a_m^{-\frac{1}{2}}$$
Fix $m=0$, and let $u^E(n)=F_{n,0}$, then according to the previous proof, $\rho-a.e.\ E$, $u^{E}$ is generalized eigenfunction of $Hu=Eu$ and obey the estimate
$$|u^E(n)|\le a_0^{-\frac{1}{2}} a_n^{-\frac{1}{2}}$$
By the choice of $a_n$, we have
$$|u^E(x_k)|\le (1+|k|)^{y}.$$ \qed
\section{Acknowledgments}
We are grateful to S. Molchanov for mentioning the extended Schnol's
Theorem in his talk at UCI: it has become an important part of our
proof. S.J. is a 2014--15 Simons Fellow. This research was partially
    supported by the NSF DMS--1401204.
    We are also grateful to the Isaac
    Newton Institute for Mathematical Sciences, Cambridge, for support
    and hospitality during the programme Periodic and Ergodic Spectral
    Problems where a part of this work was done.

\end{document}